\documentclass[numbook,runningheads]{svjour3}
\usepackage{hyperref}
\usepackage{amssymb}
\usepackage{amsmath}
\usepackage{algorithm}
\usepackage{algpseudocode}
\usepackage{comment}

\usepackage{graphicx}

\providecommand{\E}[1]{{\ensuremath{\mathrm{E}}\mspace{-2mu}\left[#1\right]}}\providecommand{\prob}[1]{{\ensuremath{\mathbb{P}}\mspace{-2mu}\left[#1\right]}}\providecommand{\var}[1]{{\ensuremath{\mathrm{Var}}\mspace{-2mu}\left[#1\right]}}\providecommand{\vg}{{\ensuremath{V_0}}}
\providecommand{\tol}{\mathrm{TOL}}
\providecommand{\rset}{\mathbb{R}}

\providecommand{\nset}{\mathbb{N}}
\providecommand{\discreteh}{\mathfrak{H}}
\providecommand{\hier}{\ensuremath{\mathbf{H}}}
\providecommand{\work}{\ensuremath{W}}

\providecommand{\Order}[1]{ {\ensuremath{ \mathcal O\left( #1 \right)}} }

\providecommand{\levsep}{ {\ensuremath{\beta }} }

\title{Optimization of mesh hierarchies in Multilevel {M}onte {C}arlo samplers}

\author{
  Abdul--Lateef~Haji--Ali
  \and Fabio~Nobile
  \and Erik~von~Schwerin
  \and Ra\'{u}l~Tempone
}

\institute{
  A.~Haji--Ali (\email{abdullateef.hajiali@kaust.edu.sa}) \and R.~Tempone  (\email{raul.tempone@kaust.edu.sa})\at
        Applied Mathematics and Computational Sciences, KAUST, Thuwal, Saudi Arabia.
  \and
  F.~Nobile  \at {MATHICSE-CSQI, EPF de Lausanne,
    Switzerland.}   \and
 E.~von~Schwerin \at {
Department of Mathematical Sciences, University of Delaware, Newark, USA}
}

\def\makeheadbox{{}}
\date{}

\begin{document}

\maketitle

\begin{abstract}
We perform a general optimization of the parameters in the
Multilevel Monte Carlo (MLMC) discretization hierarchy based on
uniform discretization methods with general approximation orders and
computational costs.
We optimize hierarchies with geometric and non-geometric sequences
of mesh sizes and show that geometric hierarchies, when optimized,
are nearly optimal and have the same asymptotic computational
complexity as non-geometric optimal hierarchies. We discuss how
enforcing constraints on parameters of MLMC hierarchies affects the
optimality of these hierarchies. These constraints include an upper
and a lower bound on the mesh size or enforcing that the number of
samples and the number of discretization elements are integers. We
also discuss the optimal tolerance splitting between the bias and
the statistical error contributions and its asymptotic behavior. To
provide numerical grounds for our theoretical results, we apply
these optimized hierarchies together with the Continuation MLMC
Algorithm \cite{haji_CMLMC}.  The first example considers a
three-dimensional elliptic partial differential equation with random
inputs.  Its space discretization is based on continuous piecewise
trilinear finite elements and the corresponding linear system is
solved by either a direct or an iterative solver.  The second
example considers a one-dimensional It\^o stochastic differential
equation discretized by a Milstein scheme.
\end{abstract}



\journalname{SPDE}

\section{Introduction}\label{s:intro}\makeatletter{}The history of Multilevel Monte Carlo methods can be traced back to
Heinrich et al.~\cite{heinrich98,hs99}, where it was introduced in the
context of parametric integration.  Kebaier~\cite{kebaier05} then used
similar ideas for a two-level Monte Carlo (MC) method to approximate
weak solutions to stochastic differential equations (SDEs) in
mathematical finance. The basic idea of the two-level MC method is to
reduce the number of samples on the fine mesh by using a control
variate that is obtained by approximating the solution on a coarser
mesh.  In~\cite{giles08}, Giles extended this idea to more than two
levels and dubbed his extension the Multilevel Monte Carlo (MLMC)
method.  Giles introduced a hierarchy of discretizations with
geometrically decreasing mesh sizes. His work also included an
optimization of the number of samples on each level that reduced the
computational complexity to $\Order{\tol^{-2} (\log{\tol})^2}$ when
applied to SDEs with Euler-Maruyama discretization, compared to
$\Order{\tol^{-3}}$ of the standard Euler-Maruyama MC method.  In
\cite{giles08_m}, Giles further reduced the computational complexity
of approximating weak solutions of a one-dimensional SDE to
$\Order{\tol^{-2}}$ by using the Milstein scheme instead of the
Euler-Maruyama scheme to discretize the SDE.  MLMC has also been
extended and applied in many contexts, including equations with jump
diffusions \cite{xg12}, partial differential equations (PDEs) with
stochastic coefficients \cite{bsz11,cst13,cgst11,tsgu13} and
stochastic partial differential equations (SPDEs) \cite{bls13,gr12},
to compute scalar quantities of interest that are functionals of the
solutions.  In \cite[Theorem 2.3]{tsgu13}, an optimal convergence rate
is derived for general rates of strong and weak convergence and the
computational complexity associated with generating a single sample of
the quantity of interest.  It is shown that if the strong convergence
is sufficiently fast, the computational complexity can be of the
optimal rate, $\Order{\tol^{-2}}$.

Several points can be investigated in this standard MLMC setting. For
instance, the standard MLMC uses uniform mesh sizes on each level and
across the levels the mesh sizes follow a geometric sequence in which
the ratio between mesh sizes of subsequent levels is a constant,
$\levsep$, henceforth referred to as level separation. However, it is
not clear if this is an optimal choice.  Moreover, in the literature,
the derivation of the optimal number of samples on each level assumed
an equal, fixed splitting of accuracy between statistical and bias
error contributions.  In \cite{haji_CMLMC}, the authors used a more
efficient splitting that improved the running time of MLMC by a
constant factor, but no analysis of the splitting parameter was
provided.  In this work, we show that, in certain cases, the optimal
level separation is not a constant and depends on several parameters,
including the level index, $\ell$.  Moreover, when restricted to
geometric, but not nested, hierarchies, we optimize for the
constant level separation parameter, $\levsep$, by using some
heuristics and show that using this optimized value the computational
complexity of the geometric hierarchies is close to the computational
complexity of the optimized non-geometric hierarchies.  We also show
that the computational complexity of both hierarchies are the same in
the limit as $\tol \to 0$.  In addition, we analyze the optimal
splitting parameter, $\theta$, and note its asymptotic behavior as
$\tol \to 0$.  Several issues arise in a practical implementation of
MLMC. One of these issues is that the hierarchies generated by
optimality theorems are usually not applicable due to constraints on
either mesh sizes (for instance due to CFL stability limitations) or
the number of samples; the constraint on the latter being an integer,
for example.  We analyze these issues and note their effect on the
optimality of the MLMC hierarchies.  Other issues include the stopping
criteria \cite{Bayer2012} and the estimation of variances in the case
of a small number of samples, a feature that is inherent to MLMC and
is always present in the deepest levels of the MLMC hierarchies.  To
this end, we here apply these optimized hierarchies together with the
Continuation MLMC algorithm (CMLMC) \cite{haji_CMLMC} and show the
effectiveness of the resulting algorithm in several examples.  The use
of {\it a posteriori} error estimates and related adaptive algorithms,
as introduced first in~\cite{hsst12}, is beyond the scope of this
work, which focuses instead on optimizing {\it a priori} defined
parametric families to create the discretization hierarchies.

This work is organized as follows.  Section~\ref{s:hier_intro} recalls
the MLMC sampling framework and states the hierarchy optimization
problem.  Several approximation steps lead to an analytically
treatable problem.  Section~\ref{s:nongeo} presents the solution for
the case of unconstrained optimal mesh sizes, including the number of
samples per level and the splitting accuracy parameter; these optimal
mesh sizes do not form geometric sequences in general.  Then,
Section~\ref{s:geo} presents the optimal hierarchies if they are
restricted to geometric sequences of mesh sizes.  Finally,
Section~\ref{s:res} illustrates the theoretical results with numerical
examples, which include three-dimensional PDEs with random inputs and
It\^o SDEs, and Section~\ref{s:conc} draws conclusions and proposes
future extensions of this work.  To avoid cluttering the presentation,
the technical derivations of the formulas included in this work are
outlined in the appendix.

\section{Optimal MLMC hierarchies}\label{s:opt}\makeatletter{}Here we state the problem of optimizing the mesh hierarchies in MLMC
and present the mesh hierarchies resulting from a theoretical
optimization, first allowing very general sequences of mesh sizes and
then for comparison restricting ourselves to geometric sequences.

In Section~\ref{s:hier_intro} we introduce the MLMC hierarchy, the
parameters that we consider free to optimize in the hierarchy, and the
models of the computational work and of the weak and strong
errors that define the general, discrete and non-convex, optimization
problem. Simplifying assumptions then lead to an analytically
treatable continuous optimization problem in Sections~\ref{s:nongeo}--\ref{s:geo}.

\subsection{Problem setting}\label{s:hier_intro}\makeatletter{}Let $g(u)$ denote a scalar quantity of interest, which is a function of
the solution $u$ of an underlying stochastic model. Our goal is to
approximate the expected value, $\E{g(u)}$, to a given accuracy $\tol$
with a high probability of success. We assume that individual outcomes
of the
underlying solution, $u$, and the evaluation of $g(u)$ are approximated
by a discretization-based numerical scheme characterized by a mesh
size\footnote{We consider uniform meshes, but the extension to
  certain non-uniform meshes is immediate; see
  Remark~\ref{rem:mesh_grading} in
  Section~\ref{s:nongeo}.}, $h$.
The following examples are adapted from \cite{haji_CMLMC} with some modification:
\begin{example}
\label{ex:spde_problem}
Let $(\Omega,\mathcal{F},\mathbb P)$ be a complete probability space
and $\mathcal{D}$ be a open, bounded and convex polygonal domain in $\rset^d$.
 Find $u: \mathcal{D} \times \Omega \to \rset$ that solves $\mathbb P$-almost surely (a.s.)
the following equation:
\begin{subequations}
\begin{align*}
    -\nabla \cdot \left( a \nabla u \right) &= f, &\text{ in }
             \mathcal D,  \\
     u &=  0, &\text{ on }  \partial \mathcal D,
\end{align*}
\end{subequations}
  where the value of the diffusion coefficient
  and the forcing are represented by random fields, yielding a random
  solution. We wish to compute $\E{g(u)}$ for some deterministic
  functional $g$ which is globally Lipschitz satisfying $|g(u) - g(v)|
  \leq G \|u-v\|_{H^1(\mathcal D)}$ for some constant $G>0$ and all
  $u,v\in H^1(\mathcal D)$.
Following \cite{tsgu13}, we also make the following assumptions
\begin{itemize}
\item $a_{\min} = \min_{\mathbf x \in \mathcal D}a(\mathbf x) > 0$
  a.s. and $1/a_{\min} \in L^p_\mathbb{P}(\Omega)$, for all
  $p \in (0, \infty)$.
  \item $a \in L^p_{\mathbb P}(\Omega, C^1(\overline{\mathcal D}))$,
    for all $p \in (0, \infty)$.
  \item $f \in L^{p^*}_{\mathbb P}(\Omega, L^2(\mathcal D))$ for some
    $p^* > 2$.
  \end{itemize}
  Here, $L^p_{\mathbb P}(\Omega, \mathcal B)$ is the Bochner space of
  $p$-th integrable $\mathcal B$-valued random fields, where the
  $p$-th integrability is with respect to the probability measure
  $\mathbb P$.  On the other hand, $C^1(\overline{\mathcal D})$ is the
  space of continuously differentiable functions with the usual norm
  \cite{cst13}.  Note that with these assumptions and since
  $\mathcal D$ is bounded, one can show that
  $\max_{\mathbf x \in \mathcal D}a(\mathbf x) < \infty$ a.s.  A
  standard approach to approximate the solution of the previous
  problem is to use finite elements on regular triangulations. In such
  a setting, the parameter $h>0$ refers to either the maximum element
  diameter or another characteristic length and the corresponding
  approximate solution is denoted by $u_h$.  For piecewise linear or
  piecewise $d$-multilinear continuous finite element approximations,
  and with the previous assumptions, it can be shown \cite[Corollary
  3.1]{tsgu13} that asymptotically as $h \to 0$:
\begin{itemize}
\item $\left|\E{g(u)-g(u_h)} \right| = \Order{h^2}$.
\item $\var{g(u)-g(u_h)} = \Order{h^4}$.
\end{itemize}

\end{example}

\begin{example}
\label{ex:sde_problem}
 Here we study the weak approximation of It\^o stochastic differential
 equations (SDEs).
Let $(\Omega,\mathcal{F},\mathbb P)$ again denote a complete probability space
 and let
\begin{equation}
  du(t) = a(t,u(t))dt +
  b(t,u(t))dB(t), \qquad 0<t<T,
  \label{eq:sde}
\end{equation}
where $u(t)$ is a stochastic process in $\rset^d$, with
randomness generated by a $k$-dimensional Wiener process with
independent components, $B(t)$, cf.~\cite{KS,Ok}, and
$a(t,u)\in\rset^d$ and $b(t,u)\in\rset^{d\times k}$ are the drift and
diffusion fluxes, respectively.  For any given sufficiently well-behaved function, $g: \rset^d\rightarrow \rset$, our goal is to approximate the expected
value, $\E{g(u(T))}$. A typical application is to compute option prices in
mathematical finance, cf.~\cite{JouCviMus_Book,Glasserman}, and other
related models based on stochastic dynamics.

When one uses a standard Milstein scheme based on uniform time steps of size $h$ to approximate \eqref{eq:sde},
the following rates of approximation hold: $E[g(u(T))-g(u_h(T))] = \Order{h}$ and
$E[(g(u(T))-g(u_h(T)))^2] = \Order{h^2}$.
For suitable assumptions on the functions $a$, $b$, and $g$, we refer to \cite{milstein2004stochastic}.

\end{example}

To avoid cluttering the notation, we omit the reference to the underlying solution from
now on, simply denoting the quantity of interest as $g$.
Following the standard MLMC approach, we introduce a hierarchy of $L+1$
meshes defined by decreasing mesh sizes $\{h_\ell\}_{\ell=0}^L$ and
we denote the resulting approximation of $g$ using mesh size $h_\ell$ by
$g_\ell$, or by $g_\ell(z)$ when we want to stress the dependence
on an outcome on some underlying random variable $z$. Then, the expected value
of the finest approximation, $g_L$,
can be expressed as
\begin{align*}
    \E{g_L} & =
  \E{g_0} + \sum_{\ell=1}^L \E{g_{\ell}-g_{\ell-1}},
\end{align*}
and the MLMC estimator is obtained by approximating the expected values
in the telescoping sum by sample averages as
\begin{align}
  \label{eq:MLMC_estimator}
  \mathcal{A} & = \frac{1}{M_0}\sum_{m=1}^{M_0}g_0(z_{0,m}) +
  \sum_{\ell=1}^L\frac{1}{M_\ell}\sum_{m=1}^{M_\ell}
  \left(g_{\ell}(z_{\ell,m})-g_{\ell-1}(z_{\ell,m})\right),
\end{align}
where, for every $\ell$,  $\{z_{\ell,m}\}_{m=1}^{M_\ell}$
denotes independent identically distributed (i.i.d.) random variables
representing the underlying, mesh-independent, stochastic model.
In addition, the random variables in the union of all these sets are independent.
We note that, given the model for $g_\ell$, the MLMC estimator is
defined by the triplet
$\hier = \left(L,\{h_\ell\}_{\ell=0}^L, \{M_\ell\}_{\ell=0}^L\right)$,
which we also refer to as the MLMC hierarchy.  Depending on the
numerical discretization method, possible mesh sizes will be
restricted to a discrete set of positive real numbers, which we denote
by $\discreteh$.  For instance, for uniform meshes in the domain
$(0,1)^d$, the number of subdivisions in each dimension has to be an
positive integer, resulting in the constraint $h^{-1} \in \nset_+$.
We do not, however, introduce any other restriction on the mesh sizes
but allow the MLMC hierarchy to use any decreasing sequence of
attainable mesh sizes.
Moreover, the number of samples on any level is a positive integer,
$M_\ell\in\nset_+$, while $L$ is a non-negative integer, $L\in\nset$.

If $W_\ell$ is the average cost associated with generating one sample of the difference,
$g_{\ell}-g_{\ell-1}$, or simply $g_0$ if $\ell=0$, then the cost of the
estimator~\eqref{eq:MLMC_estimator} is
\begin{align}
  \label{eq:work_sum}
  \work(\hier) & = \sum_{\ell=0}^{L} M_\ell W_\ell.
\end{align}
We assume that the work required to generate one sample of mesh size
$h$ is proportional to $h^{-d\gamma}$, where $d$ is the dimension of the
computational domain and $\gamma > 0$ represents the complexity of
generating one sample with respect to the number of degrees of freedom.
Thus, we model the average cost on level $\ell$ as
\begin{align}
  \label{eq:wl_model}
  W_\ell &= h_\ell^{-d\gamma},
\end{align}
and consequently use the representation
\begin{align}
  \label{eq:model_work}
  \work(\hier) & = \sum_{\ell=0}^{L} \frac{M_\ell}{h_\ell^{d\gamma}}
\end{align}
for the total work to evaluate the MLMC estimator~\eqref{eq:MLMC_estimator}.
This can be motivated in two ways. Namely, we are simply neglecting the work to generate the coarser variable in each
realization pair $(g_\ell,g_{\ell-1})$ or, we are bounding the work to generate the pair by a constant factor, which is clearly
less than or equal to twice the work to generate the finest variable in each realization pair.

For example, if each sample evaluation is the approximation of an
It\^{o} stochastic differential equation by a time stepping scheme,
then $d=\gamma=1$. If, instead, the underlying differential
equation is an elliptic partial differential equation with a
stochastic coefficient field, then a numerical method based on an ideal
multigrid solver will still have $\gamma=1$ up to a logarithmic
factor, while a naive implementation of Gaussian elimination based on full matrices leads to
$\gamma=3$.

We want to find a hierarchy, $\hier$, which, with a prescribed failure
probability, $0 < \alpha \ll 1$, satisfies
\begin{align*}
    \prob{\left|\E{g}-\mathcal{A}\right| > \tol} \leq \alpha,
\end{align*}
while minimizing the work, $\work(\hier)$.
Here, we aim to meet this accuracy requirement by
controlling the bias and statistical error separately as
\begin{align}
  \label{eq:splitting}
  \left|\E{g-\mathcal{A}}\right| \leq (1-\theta)\tol \quad
  & \quad\quad\text{and} &
  \prob{\left|\E{\mathcal{A}}-\mathcal{A}\right| > \theta\,\tol} \leq
                           \alpha.
\end{align}
This splitting of the error introduces a new parameter, $0<\theta<1$,
which we are free to choose. We will later see that the choice of
$\theta$ that minimizes the work is not obvious, and does not
reduce to any simple rule of thumb.
Motivated by the Lindeberg-Feller Central Limit
Theorem in the limit $\tol\to 0$; see~\cite[Lemma~A.2]{haji_CMLMC},
the probabilistic constraint in \eqref{eq:splitting} can be replaced
by a constraint on the variance of the estimator as follows:
\begin{align}
  \label{eq:var_bound}
  \var{\mathcal{A}} & \leq \left(\frac{\theta \tol}{C_\alpha} \right)^2,
\end{align}
where $C_\alpha$ satisfies $\Phi(C_\alpha) = 1-\frac{\alpha}{2}$;
here, $\Phi$ is the standard normal cumulative distribution function.

By construction of the estimator, $\E{\mathcal{A}}=\E{g_L}$ and using the notation
\[
V_\ell = \begin{cases} \var{g_0} & \ell=0,\\
 \var{g_\ell - g_{\ell-1}} & \ell>0,
    \end{cases}
\]
and by independence we have $\var{\mathcal{A}} = \sum_{\ell=0}^L V_\ell M_\ell^{-1}$.
The requirements~\eqref{eq:splitting} and \eqref{eq:var_bound} therefore become
\begin{subequations}
  \label{eq:requirements}
  \begin{align}
    \label{eq:bias_requirement}
    \left|\E{g-g_L}\right| & \leq (1-\theta)\tol, \\
    \label{eq:var_requirement}
    \sum_{\ell=0}^L V_\ell M_\ell^{-1}
    & \leq
    \left(\frac{\theta \tol}{C_\alpha} \right)^2.
  \end{align}
\end{subequations}

We now assume that the numerical approximation of $g_\ell$ leads to
weak convergence of order $q_1$ and strong convergence of order
$q_2\leq 2q_1$ as $h\to 0$, and we further assume that the variance
on the coarsest level is approximately independent of its
corresponding mesh size. Note that the condition that is usually
assumed for MLMC is $\min(q_2,d\gamma) \leq 2q_1$, c.f. \cite[Theorem
2.3]{tsgu13}, to ensure that the cost of MLMC is not dominated by the
cost of a single sample on each level. We assume instead the slightly
more restrictive condition $q_2 \leq 2 q_1$, since it does not depend
on the dimensionality of the problem, $d$.
Using these assumptions and neglecting all higher order terms in
$h_\ell$, we postulate for some constants, $0<Q_W, Q_S<\infty$ the
following models for the bias and variances:
\begin{align}
  \nonumber
  \left|\E{g-g_L}\right| &= Q_W h_L^{q_1}, \\
  \label{eq:var_model}
  V_\ell  &= Q_S h_{\ell-1}^{q_2}, \quad \text{ for } \ell > 0.
\end{align}
These models are reasonable for many cases (such as those listed in
Section~\ref{s:res}) and hold asymptotically as $\ell \to \infty$. We
limit ourselves to these cases,
while acknowledging that there are others, such as those in
\cite{moraes2014multilevel,tesei2014multi}, that do not follow the model \eqref{eq:var_model}.
 We observe that the problem of finding
$\hier=(L,\{h_\ell\}_{\ell=0}^L, \{M_\ell\}_{\ell=0}^L)
\in\nset\times\discreteh^{L+1}\times\nset_+^{L+1}$
minimizing $\work(\hier)$ in~\eqref{eq:model_work} while satisfying
the constraints~\eqref{eq:requirements} is a difficult discrete
optimization problem.  Hence, we make a further simplification by
temporarily removing the constraints on $h_\ell$ and $M_\ell$ to let
$\hier\in\nset\times\rset_+^{L+1}\times\rset_+^{L+1}$. The simplified
variance model~\eqref{eq:var_model} is valid for example for nested
geometric sequences of mesh sizes, but in our more general setting in
this paper it can instead be seen as a penalty on closely spaced
meshes, where it overestimates the resulting variance.

The simplified models for the bias and the variance of the MLMC
estimator are then
\begin{subequations}
  \label{eq:model_error}
  \begin{align}
    \label{eq:model_bias}
    \left|\E{g-\mathcal{A}}\right| &= Q_W h_L^{q_1}, \\
    \label{eq:model_variance}
    \var{\mathcal{A}} &= \frac{\vg}{M_0} +
    Q_S\sum_{\ell=1}^L\frac{h_{\ell-1}^{q_2}}{M_\ell},
  \end{align}
\end{subequations}
with problem- and method-specific positive constants, $Q_W$, $Q_S$, and $V_0$.
We note that neglecting the higher-order terms in $h_\ell$ is usually
justified in the model of the bias, which only depends on the finest
mesh. On the other hand, it may be that the contribution of the higher-order terms on
coarse meshes makes the model \eqref{eq:model_variance} inaccurate,
       causing the hierarchies derived in this work to be suboptimal.
    Dealing with such non-asymptotic behavior is beyond the scope of this work and we leave it for future work.

\subsection{General mesh size sequences}\label{s:nongeo}\makeatletter{}    Here, we present the optimal hierarchy, $\hier$, using the
continuous, convex, model of the previous subsection, which solves:
\begin{problem}  \label{prob:cont_opt}
  Find $\hier = (L,\{h_\ell\}_{\ell=0}^L,\{M_\ell\}_{\ell=0}^L)
  \in \nset\times\rset_+^{L+1}\times\rset_+^{L+1}$ such that
  \begin{subequations}
    \label{eq:cont_opt}
    \begin{align}
      \label{eq:objective}
      \work(\hier) & =
      \sum_{\ell=0}^{L} \frac{M_\ell}{h_\ell^{d\gamma}},
    \end{align}
        is minimized while satisfying the constraints    \begin{align}
      \label{eq:bias_constraint}
      Q_W h_L^{q_1} & \leq (1-\theta) \tol, \\
      \label{eq:variance_constraint}
      \frac{\vg}{M_0} +
      Q_S\sum_{\ell=1}^L\frac{h_{\ell-1}^{q_2}}{M_\ell}
      & \leq \left(\frac{\theta \tol}{C_\alpha} \right)^2,
    \end{align}
  \end{subequations}
  for some $\theta\in(0,1)$.
\end{problem}

Note that, even though the parameter $\theta$ is not part of the
hierarchy $\hier$ defining the MLMC estimator, determining
$\theta$ is still an important part of the optimization.
Initially, we treat the parameters $\theta$ and $L$ as given and
optimize first with respect to $\{M_\ell\}_{\ell=0}^L$ and then
$\{h_\ell\}_{\ell=0}^L$.
From a Lagrangian formulation of the problem of minimizing the general
work model~\eqref{eq:work_sum} under the
constraint~\eqref{eq:var_requirement}, it is straightforward to obtain
the optimal number of samples,
\begin{align}
  \label{eq:opt_M_ell}
  M_\ell & = \left( \frac{C_\alpha}{ \theta\tol} \right)^2
  \sqrt{\frac{V_\ell}{W_\ell}} \sum_{k=0}^L \sqrt{W_k V_k},
\end{align}
in terms of general work estimates, $\{W_\ell\}_{\ell=0}^L$, and
variance estimates, $\{V_\ell\}_{\ell=0}^L$; see Section~\ref{s:app_opt_hier}
for more details on this and the following steps.
The finest mesh size is determined by the bias
constraint~\eqref{eq:bias_constraint}, for any given choice of~$\theta$.
The optimality conditions then lead to a linear difference equation, which
can easily be solved for the remaining mesh sizes.
In the idealized situation, where the coarsest mesh size is treated as
an unconstrained variable in the optimization,
we can analytically minimize the computational complexity with respect to $\theta$ to obtain the
optimal hierarchy for any fixed $L$.
Introducing the two model- and method-dependent parameters,
\begin{align}
  \label{eq:def_chi_eta}
  \eta  = \frac{q_1}{d\gamma} & \qquad\qquad\text{and} & \chi = \frac{q_2}{d\gamma},
\end{align}
we can summarize the result derived in Sections~\ref{s:app_chi_1} and~\ref{s:app_chi_n1} in the following theorems for the two cases: $\chi=1$ and $\chi \neq 1$.
\begin{theorem}[On the optimal hierarchies when $\chi=1$]
  \label{thm:opt_hier_chi_1}
    For any fixed $L\in\nset$, with $\chi = 1$, the optimal
  sequences $\{h_\ell\}_{\ell=0}^L$ and $\{M_\ell\}_{\ell=0}^L$ in
  Problem~\ref{prob:cont_opt} are given by
\begin{subequations}
  \begin{align}
    \label{eq:hl_optimal_chi1}
    h_\ell & = \levsep^{L-\ell}
    \left(\frac{(1-\theta)\tol}{Q_W}\right)^{\frac{1}{q_1}},
    && \text{for $l =0,1,2,\ldots,L$,} \\
    \label{eq:Ml_optimal_chi1}
    M_\ell & =
    \levsep^{-q_2\ell}\vg(L+1)\left(\frac{C_\alpha}{\theta\tol}\right)^2,
    && \text{for $l =0,1,2,\ldots,L$,}
    \intertext{where the level separation $\levsep\in(1,\infty)$ is
      independent of $\ell$,}
    \label{eq:opt_beta_chi1}
    \levsep & = \left\{
      \left(\frac{Q_W}{(1-\theta)\tol}\right)^\frac{1}{q_1}
      \left(\frac{\vg}{Q_S}\right)^\frac{1}{q_2}
    \right\}^{\frac{1}{L+1}}, &&
  \end{align}
  and the optimal choice of the splitting parameter
  \begin{align}
  \label{eq:opt_theta_chi1}
  \theta(1, \eta, L) & =
  \left(
    1 + \frac{1}{2\eta}\frac{1}{L+1}
  \right)^{-1}.
  \end{align}
\end{subequations}
\end{theorem}
\begin{lemma}
  \label{thm:opt_L_chi_1}
  For the case $\chi=1$ and the optimal hierarchies in
  Theorem~\ref{thm:opt_hier_chi_1}, the optimal number of levels, $L$,
  satisfies
  \begin{align}
    \label{eq:L_chi_1_bounds} 1  <
    \frac{2\eta(L+1)}
    {\log\left(\tol^{-1} Q_W V_0^\eta Q_S^{-\eta} \right)} \leq
    \frac{1}{\exp(1)-1}
    \left(
      \exp(1) + \frac{1}{\log{\left(\tol^{-1} Q_W V_0^\eta Q_S^{-\eta}\right)}}
    \right),
    \end{align}
  for any $\tol<Q_W V_0^\eta Q_S^{-\eta}$, and asymptotically
  \begin{align}
        \label{eq:L_chi_1_asymb}  \lim_{\tol \to 0} \frac{L+1}{\log{\tol^{-1}}} = \frac{1}{2\eta}.
    \end{align}
\end{lemma}
\begin{corollary}
  \label{thm:opt_work_chi_1}
  For the case $\chi=1$ and the optimal hierarchies in Theorem~\ref{thm:opt_hier_chi_1} and using $L$ in \eqref{eq:L_chi_1_asymb}, the total work \eqref{eq:model_work} satisfies
    \begin{equation}
  \label{eq:opt_asymp_work_chi1}
  \frac{W(\hier)}{\tol^{-2}(\log{\tol})^2} \to  C_\alpha^2 \exp(2) Q_S \left(\frac{1}{2\eta}\right)^2,\quad \text{ as } \tol \searrow 0.
      \end{equation}
\end{corollary}

\begin{theorem}[On the optimal hierarchies when $\chi\neq 1$]
    \label{thm:opt_hier_chi_n1}
    For any fixed $L\in\nset$, with $\chi\neq 1$, the optimal
    sequences, $\{h_\ell\}_{\ell=0}^L$ and $\{M_\ell\}_{\ell=0}^L$, in
    Problem~\ref{prob:cont_opt} are given by
    \begin{subequations}
        \label{eq:parameters}
        \begin{equation}
        \label{eq:hl_optimal}
        \begin{split}
            h_\ell(\theta, L) = &
            \left(
            \frac{(1-\theta)\,\tol}{Q_W}
            \right)^{\frac{1}{q_1}\frac{1-\chi^{\ell+1}}{1-\chi^{L+1}}}
            \left(
            \frac{\vg}{Q_S}
            \right)^{\frac{1}{d\gamma}\frac{\chi^\ell-\chi^L}{1-\chi^{L+1}}}
            \quad\quad\quad\quad\quad\quad\,\,\\
            &\cdot
            \chi^{-\frac{1}{d\gamma}\frac{2}{1-\chi}
            \left(
                \frac{\chi^{L+1}-\chi^{\ell+1}}{1-\chi^{L+1}} +
                \frac{L\left(1-\chi^{\ell+1}\right)
                -\ell\left(1-\chi^{L+1}\right)}
                {1-\chi^{L+1}}
            \right)
            },                    \end{split}
    \end{equation}
    \begin{equation}
        \label{eq:Ml_optimal}
        \begin{split}
        M_\ell(\theta, L) = &
        \left(\frac{C_\alpha}{\theta\tol}\right)^2
        \left((1-\theta)\,\tol\right)^{\frac{\chi}{\eta}\frac{1-\chi^\ell}{1-\chi^{L+1}}}
        \vg^{\frac{\chi^\ell-\chi^{L+1}}{1-\chi^{L+1}}}\quad\quad\quad\quad\\
        &\cdot
        \left(
        \frac{Q_S^{1/\chi}}
        {Q_W^{1/\eta}}
        \right)^{\frac{\chi(1-\chi^\ell)}{1-\chi^{L+1}}}
        \frac{1-\chi^{L+1}}{\chi^L(1-\chi)}
        \chi^{\left\{-\frac{2\chi}{1-\chi}\frac{1-\chi^\ell}{1-\chi^{L+1}}(L+1)
            + \frac{1+\chi}{1-\chi}\ell\right\}},        \end{split}
    \end{equation}
    where the optimal choice of the splitting parameter is
    \begin{equation}
            \label{eq:theta_opt}
        \theta(\chi, \eta, L) =
        \left(
            1 + \frac{1}{2\eta}\frac{1-\chi}{1-\chi^{L+1}}
        \right)^{-1}.
    \end{equation}
    \end{subequations}
\end{theorem}
\begin{lemma}
    \label{thm:opt_L_chi_n1}
    For the case $\chi\neq 1$ and the optimal hierarchies in Theorem~\ref{thm:opt_hier_chi_n1}, the optimal number of levels, $L$, satisfies
    \begin{align}
    \label{eq:L_bound_TOL}
    \frac{1}{c_2}
    \left(1+\frac{c_1+\log{\left(1+2\eta\right)}}{\log{\left(\tol^{-1}\right)}}\right)
    \! < \! \frac{L+1}{\log{\left(\tol^{-1}\right)}}
    & \! < \!
    \begin{cases}
        \frac{1}{c_2}
        \left(1+\frac{c_1+\log{\left(1+\frac{2\eta}{1-\chi}\right)}}
        {\log{\left(\tol^{-1}\right)}}\right), & \!\!\!\!\!\!\chi\in(0,1),
        \\
        \frac{\chi}{c_2}
        \left(1+\frac{c_1+\log{\left(\frac{2\eta}{\chi-1}\right)}}
        {\log{\left(\tol^{-1}\right)}}\right), & \!\!\!\!\!\!\chi\in(1,\infty),
    \end{cases}
    \end{align}
    where
    \begin{align}
    \label{eq:L_bounds_consts}
    c_1 = \log{\left(\frac{\vg^{\eta/\chi}}{Q_S^{\eta/\chi}}Q_W\right)} \qquad\text{and}\qquad
    c_2 = \log{(\chi)}\frac{2\eta}{\chi-1} > 0,
    \end{align}
    and asymptotically
  \begin{align}
    \label{eq:L_bounds}
    \frac{1}{2\eta}\frac{\chi-1}{\log{\chi}}
    \leq \liminf_{\tol\to 0}\frac{L+1}{\log{\left(\tol^{-1}\right)}}
    & \leq \limsup_{\tol\to 0}\frac{L+1}{\log{\left(\tol^{-1}\right)}}
    \leq
    \frac{\max{\{1,\chi\}}}{2\eta}\frac{\chi-1}{\log{\chi}}.
  \end{align}
\end{lemma}

\begin{corollary}
    \label{thm:opt_work_chi_n1}
    For the case $\chi\neq 1$ and the optimal hierarchies in Theorem~\ref{thm:opt_hier_chi_n1} and using the upper bound on $L$ in \eqref{eq:L_bound_TOL}, the total work \eqref{eq:model_work} satisfies
    \begin{subequations}
        \label{eq:optimal_work}
        \begin{align}
        \label{eq:optimal_work_1}
        \frac{\work(\hier)}{\tol^{-2\left(1+\frac{1-\chi}{2\eta}\right)}}
        & \to C_1, &&
        \text{as $\tol\searrow 0$ for $\chi\in(0,1)$, and} \\
        \label{eq:optimal_work_2}
        \frac{\work(\hier)}{\tol^{-2}}
        & \to C_2, &&
        \text{as $\tol\searrow 0$ for $\chi>1$,}
        \end{align}
    \end{subequations}
    with known constants of proportionality,
    \begin{subequations}
        \label{eq:optimal_work_constants}
        \begin{align}
        \label{eq:optimal_work_C1}
        C_1 & =
        C_\alpha^2 \, Q_S Q_W^{\left\{\frac{1-\chi}{\eta}\right\}}
        \chi^{\left\{-\frac{2\chi}{1-\chi}\right\}}
        \left(\frac{1}{2\eta}\right)^2
        \left(
            1+\frac{2\eta}{1-\chi}
        \right)^{2\left(1+\frac{1-\chi}{2\eta}\right)},
        \\
        \label{eq:optimal_work_C2}
        C_2 & =
        C_\alpha^2 \, \vg^{\left\{\frac{\chi-1}{\chi}\right\}}
        Q_S^{\left\{\frac{1}{\chi}\right\}}
        \chi^{2\left\{\frac{\chi}{\chi-1}\right\}}
        \left(\chi-1\right)^{-2}.
        \end{align}
    \end{subequations}
\end{corollary}

Note that the parameter $\theta$ controlling the split between the
statistical and discretization errors depends non-trivially on the problem
parameters.
The above theorem shows that the choice of $\theta=1/2$,
used for example in the initial works by Giles~\cite{giles08,giles08_m} and
by some of the authors of the present work in~\cite{hsst13,hsst12} for adaptive
MLMC, is increasingly suboptimal as the number of levels increases.
To further understand the splitting parameter, $\theta$, we consider the asymptotic behavior
as $L(\tol)\to\infty$ and see that
\begin{subequations}
    \label{eq:theta_opt_A}
  \begin{align}
    \label{eq:theta_opt_A_I}
    \theta(\chi, \eta, L) & \to 1, && \text{as $L\to\infty$, if $\chi \geq 1$,} \\
    \label{eq:theta_opt_A_II}
    \theta(\chi, \eta, L) & \to \frac{1}{1+\frac{1-\chi}{2\eta}},
    && \text{as $L\to\infty$, if $\chi < 1$.}
  \end{align}
\end{subequations}
The qualitative observations here are: 1) if the strong
convergence is sufficiently fast, that is $\chi \geq 1$, almost all the
tolerance is allocated to the statistical error (forcing the
discretization to be fine), and 2) for slower strong convergence,
$\chi < 1$, the tolerance can be shifted either towards the
statistical error or towards the bias according to
\begin{align*}
  \lim_{L\to\infty} \theta(\chi, \eta, L) & > \frac{1}{2}\text{ (stat. error larger) },
  && \text{if $\chi < 1 < \chi+2\eta$,}\\
  \lim_{L\to\infty} \theta(\chi, \eta, L) & < \frac{1}{2}\text{ (stat. error smaller) },
  && \text{if $\chi+2\eta < 1$.}
\end{align*}
Finally, we note that the value of the optimal splitting parameter, $\theta$, in \eqref{eq:opt_theta_chi1} for $L=0$ is consistent with the single level adaptive Monte Carlo analysis in~\cite{Moon_05_II}.

Since the above theorems give the optimal $\{h_\ell\}_{\ell=0}^L$ and
$\{M_\ell\}_{\ell=0}^L$ for any given $L\in\nset$, it is easy to
find the optimal $L$ by doing an extensive search over a finite
range of integer values.  In typical cases, for computationally feasible tolerances, $L$ is a
small non-negative integer, $0\leq L\leq 10$; we can also use the obtained bounds on the optimal value of $L$ to delimit the range of possible integer values.
Moreover, using the optimal sequences $\{h_\ell\}_{\ell=0}^L$ and $\{M_\ell\}_{\ell=0}^L$ for any given $L$, we have observed that the total computational complexity is usually
rather insensitive to the value of $L$ near the optimum.

We observe that the rates in the asymptotic complexity in
Corollaries~\ref{thm:opt_work_chi_1} and~\ref{thm:opt_work_chi_n1} are the same ones obtained with more
restrictive assumptions on the sequences of mesh sizes; see for instance \cite[Theorem~1]{cgst11} and Section~\ref{s:geo}.
With the optimal number of levels, the optimized hierarchies minimize the multiplicative constants in the complexity without improving the rate.
In Corollary~\ref{thm:opt_work_chi_n1}, the blow up of the constants $C_1$ and $C_2$ as $\chi\to1$
corresponds to the need for including the
$\log(\tol^{-1})^2$ factor that appears in
the complexity of MLMC when $\chi=1$ as Corollary~\ref{thm:opt_work_chi_1} shows.

The ratio between two successive mesh sizes in Theorem~\ref{thm:opt_hier_chi_n1} has the following complicated, non-constant expression:
\begin{equation}
  \label{eq:h_ratio}
  \frac{h_{\ell+1}}{h_\ell } =
  \left(
    \frac{\vg}{Q_S}
  \right)^{-\frac{(\chi-1)\chi^\ell}{d\gamma(\chi^{L+1}-1)}}
  \chi^{\frac{2}{d \gamma}
    \left( \frac{1}{1-\chi} +
      \frac{(L+1)\chi^{\ell+1}}{\chi^{L+1}-1}\right)}
    \left( \frac{(1-\theta)\tol}{Q_W} \right)^{\frac{(\chi-1)
        \chi^{\ell+1}}{q_1 (\chi^{L+1}-1)}}.
\end{equation}

\begin{remark}[On relation to geometric hierarchies]
\label{rem:geo_hier}
  Clearly, when $\chi\neq 1$, the optimal mesh sequences are not
  geometric in general. On the other hand, according
  to~Theorem~\ref{thm:opt_hier_chi_1}, the optimal mesh sequences are
  indeed geometric when $\chi=1$. We further note that
 an asymptotic analysis when $\tol\to0$, using optimal
    $\theta$ and $L$, shows that for $L$ sufficiently large,
  the ratios~\eqref{eq:h_ratio} are
  approximately the constant $\chi^{\frac{2}{d\gamma(\chi-1)}}$
  over most of the range of $\ell$-values.
  In case $\chi<1$, this holds for  $\ell\gg 0$, and in case
  $\chi > 1$, for $\ell \ll L$. In both cases, these are the levels where
  most of the computational work would be spent using
    geometrically spaced meshes. This suggests that a
  geometric hierarchy with this level separation constant can be
  nearly optimal. We show that this is the case in
  Section~\ref{s:geo}.

\end{remark}

\begin{remark}[On non-uniform meshes]\label{rem:mesh_grading}
  The optimization and the resulting optimal hierarchies do not depend
  on the assumption that the discretizations were uniform. Indeed,
  $h_\ell$ can also be interpreted as a more general mesh parameter
  that defines a mesh size, $\Delta x_\ell$, of the underlying
  discretization as
  \begin{displaymath}
    \Delta x_\ell =  r(h_\ell, x),
  \end{displaymath}
  for some mesh grading function $r(h_\ell,x)$, allowing for example,
  for local a priori refinement of meshes close to known singularities
  in the computational domain.
  As long as approximate models~\eqref{eq:model_work}
  and~\eqref{eq:model_error} can be provided in terms of the mesh
  parameter, the expressions for the optimal hierarchies in
  Theorems~\ref{thm:opt_hier_chi_1} and~\ref{thm:opt_hier_chi_n1} can still be applied.
   As mentioned previously,
   the construction of MLMC hierarchies based on the use of a posteriori error estimates and related adaptive algorithms, as introduced first in~\cite{hsst12}, is out of the scope of the present work.
\end{remark}

\begin{remark}[On a lower bound on possible mesh sizes]
  \label{rem:h_constrained_below}
  Since equations~\eqref{eq:hl_optimal_chi1}--\eqref{eq:opt_beta_chi1}
  and~\eqref{eq:hl_optimal}--\eqref{eq:Ml_optimal}
  are
  expressed in terms of general $\theta$ and $L$, they remain valid
  when an additional constraint is imposed on the smallest
  possible mesh sizes. If for example the available computer memory dictates a
  lower limit on the practical mesh sizes, $h_\ell\geq h_\mathrm{min}$,
  then the optimal splitting for given $L$ is
  \begin{equation}
    \label{eq:theta_opt_2}
    \theta(\chi, \eta, L) =
    \begin{cases}
    \min\left\{
      1-\frac{Q_W h_\mathrm{min}^{q_1}}{\tol},
    \left(
        1 + \frac{1}{2\eta}\frac{1-\chi}{1-\chi^{L+1}}
    \right)^{-1}
    \right\},  & \text{if }\chi \neq 1,\\
    \min\left\{
      1-\frac{Q_W h_\mathrm{min}^{q_1}}{\tol},
    \left(
        1 + \frac{1}{2\eta(L+1)}
    \right)^{-1}
    \right\},  & \text{if }\chi =1,
    \end{cases}
  \end{equation}
  where tolerances $\tol\leq Q_W h_\mathrm{min}^{q_1}$ are out of
  reach of the computation. Such an extra constraint can in turn cause
  the optimal number of levels to be smaller than the lower bound
  in~\eqref{eq:L_bound_TOL} or ~\eqref{eq:L_chi_1_bounds},
  but it can still easily be found by an extensive search over a small integer set;
    the asymptotic bounds~\eqref{eq:L_chi_1_asymb} and \eqref{eq:L_bounds} are obviously not relevant then.
\end{remark}

\begin{remark}[On an upper bound on possible mesh sizes]
  \label{rem:h_constrained_above}
    If the coarsest meshes in~\eqref{eq:hl_optimal} or \eqref{eq:hl_optimal_chi1} are unfeasibly
  large for the given method of discretization, for instance due to CFL stability constraints,
  or the asymptotic models that we assumed are only valid for small enough $h_0$,
  then we should treat the largest mesh size, $h_0$, as fixed.
  We briefly analyze this case at the end of Section~\ref{s:app_chi_n1} for the case $\chi \neq 1$.
  There, we can still express all remaining mesh sizes in terms of $h_0$ and $h_L$
  by~\eqref{eq:h_inner}, and use~\eqref{eq:opt_M_ell} for the optimal
  number of samples on the resulting sequence of mesh sizes. However,
  we no longer have an explicit expression for the optimal splitting
  parameter, but only bounds from below and above
  in~\eqref{eq:double_bound_theta}. Since $L$ varies over a finite
  integer range, we can easily obtain the optimal $\theta$ and $L$ in
  a two-stage numerical optimization.
\end{remark}

\begin{remark}
  \label{rem:int_constraints_opt}
    The optimized $h_\ell$ in \eqref{eq:hl_optimal_chi1} and \eqref{eq:hl_optimal} do not necessarily belong to $\discreteh$ and might be unusable in an actual computation.
    We instead use the closest element in $\discreteh$ to each $h_\ell$. For example, for uniform meshes in the domain $(0,1)^d$ where $h_\ell^{-1}$ is the number of elements along every dimension,
    we can simply round $h_\ell^{-1}$ up to the nearest integer.
    Similarly, $M_\ell$ in \eqref{eq:Ml_optimal_chi1} and \eqref{eq:Ml_optimal} or equivalently \eqref{eq:opt_M_ell} is not necessarily an integer and we round these expression up to the nearest integer
    to get an integer number of samples that can be used in actual computations;
    see also Remark~\ref{rem:int_constraints_geo}.
  \end{remark}

\subsection{Geometric mesh size sequences}\label{s:geo}\makeatletter{}In the optimal hierarchies of Problem~\ref{prob:cont_opt} presented
above, the mesh sizes do not form a geometric sequence except for the case
$\chi=1$. In this section, we optimize MLMC hierarchies with the more restrictive
assumption that the mesh sizes \emph{do} form a geometric sequence; that is,
$h_\ell=h_0\levsep^{-\ell}$ for some positive value $\levsep>1$ and a
given $h_0$.
The work and variance models in this case become
\begin{subequations}
\label{eq:work_var_geo_models}
\begin{align}
    V_\ell &= \begin{cases}
    V_0 & \ell = 0, \\
    Q_S h_0^{q_2} \levsep^{q_2} \levsep^{-q_2\ell} & \ell > 0,
    \end{cases} \\
    W_\ell &= h_0^{-d\gamma}\levsep^{d\gamma\ell}.
\end{align}
\end{subequations}
We do \emph{not} force $\levsep$ to be a positive integer
corresponding to successive refinements of existing meshes
but instead propose
the following value of $\levsep\in(1,\infty)$:
\begin{align}
  \label{eq:opt_beta}
  \levsep &= \begin{cases}
    \chi^{\frac{2}{d\gamma(\chi-1)}}, & \text{if } \chi\in\rset_{+}\setminus\{1\},\\
    \exp\left(\frac{2}{q_2}\right), & \text{if } \chi = 1,
    \end{cases}
\end{align}
We get this value using the asymptotic analysis in Remark~\ref{rem:geo_hier} or
a heuristic optimization that treats $L$ as a real parameter (cf. Section~\ref{ss:app_geo}).
The following corollary shows the asymptotic computational complexity
of these geometric hierarchies.
\begin{corollary}
  \label{thm:opt_geo_work}
  Consider geometric hierarchies, $h_\ell=h_0\levsep^{-\ell}$, for a given $h_0$, and the optimal number of samples $M_\ell$ in
  \eqref{eq:opt_M_ell} and the work and variance models \eqref{eq:work_var_geo_models}. Moreover, assume that
  we choose $\levsep$ in \eqref{eq:opt_beta} and choose the number of
  levels, $L$, to be the lower bound of \eqref{eq:L_geo}. In other
  words, choose
\begin{align}
  \label{eq:L_geo_ceil}
  L &= \left\lceil\frac{\log{(h_0)} - \frac{1}{q_1}\log{\left(\frac{(1-\theta)\tol}{Q_W}\right)}
      }{\log(\levsep)} \right\rceil.
\end{align}
We distinguish between two cases:
  \begin{itemize}
    \item If $\chi = 1$, the optimal $\theta$ goes to $1$ as $L \to \infty,$ and the total work satisfies \eqref{eq:opt_asymp_work_chi1}.
  \item Otherwise, if $\chi \neq 1$, the optimal $\theta$ satisfies \eqref{eq:theta_opt_A} and the total work satisfies \eqref{eq:optimal_work}
    with $C_1$ as defined in \eqref{eq:optimal_work_C1} and
    \begin{align}
        \label{eq:C_2_geo}
            C_2 &= C_\alpha^2 h_0^{d\gamma(\chi-1)} \left(\sqrt{V_0} h_0^\frac{-q_2}{2} +  \sqrt{Q_S} \frac{\chi^\frac{\chi}{\chi-1}}{\chi-1} \right)^2.
    \end{align}
    Moreover, if we choose
    \begin{equation} \label{eq:geo_opt_h0} h_0 =
      \left(\frac{V_0}{Q_S}\right)^\frac{1}{q_2}
      \chi^\frac{2}{d\gamma(1-\chi)},
    \end{equation}
    then $C_2$ simplifies to \eqref{eq:optimal_work_C2}.
    Notice that~\eqref{eq:geo_opt_h0} is the limit of
      $h_0$ in Theorem~\ref{thm:opt_hier_chi_n1} when $\tol\to0$.
  \end{itemize}
\end{corollary}

\begin{remark}
\label{rem:asymb_work}
    Corollary~\ref{thm:opt_geo_work} shows that, asymptotically as $\tol \to 0$,
    the work and optimal splitting of the geometric hierarchies with optimal $\levsep$ \eqref{eq:opt_beta} is exactly the same as the work and optimal splitting of the optimized hierarchies as
    stated in Corollaries~\ref{thm:opt_work_chi_1} and~\ref{thm:opt_work_chi_n1}.
\end{remark}

\begin{figure}
\centering
\includegraphics[scale=0.5]{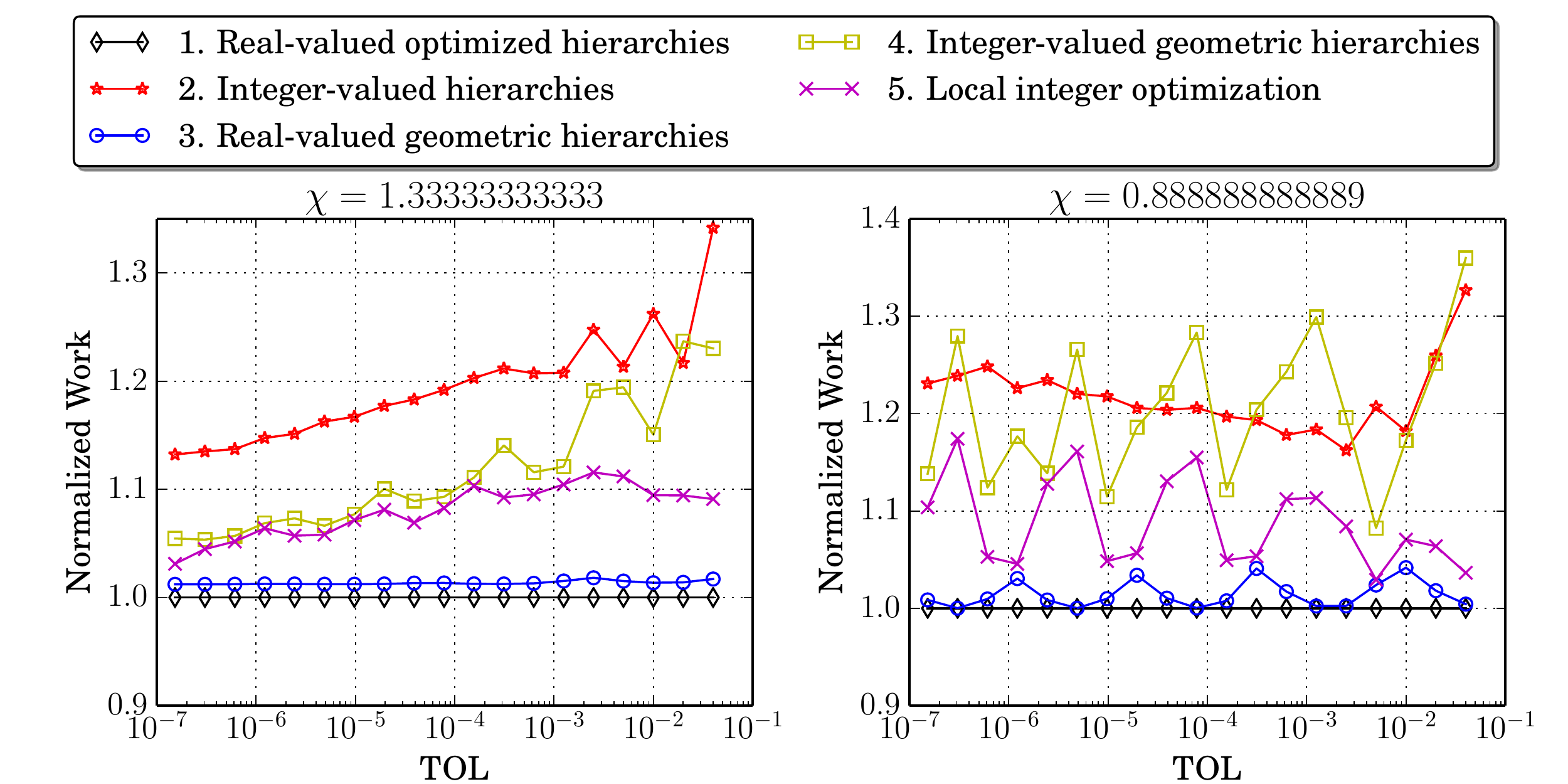}
\caption{\textbf{Ex.1}: Work, according to~\eqref{eq:model_work}, of
  different hierarchies normalized by the work estimate of the
  ``real-valued'' optimized hierarchy.  Taking the ceiling of
  $h_\ell^{-1}$ and $M_\ell$ seems to produce near-optimal
  hierarchies.  To generate these hierarchies we used the parameters
  of \textbf{Ex.1} in Table~\ref{tbl:problem_params} (See
  Remark~\ref{rem:int_constraints_geo}).  }
\label{fig:int_work}
\end{figure}

\begin{figure}
\centering
\includegraphics[scale=0.5]{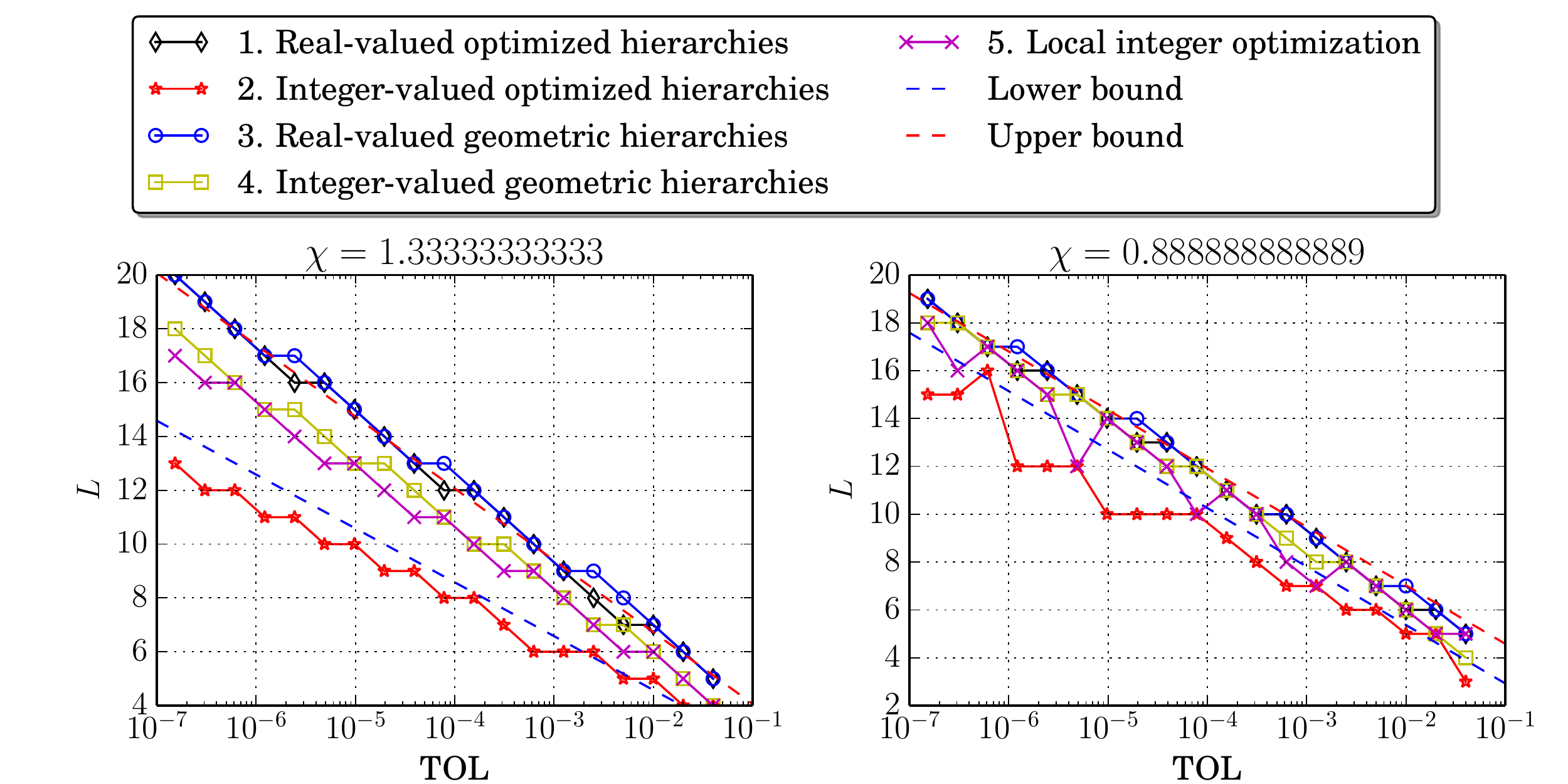}
\caption{\textbf{Ex.1}: Optimal $L$ of different hierarchies.  Here
  the bounds are from \eqref{eq:L_bound_TOL}.  To generate these
  hierarchies we used the parameters of \textbf{Ex.1} in
  Table~\ref{tbl:problem_params} (See
  Remark~\ref{rem:int_constraints_geo}).  }
\label{fig:int_L}
\end{figure}

\begin{remark}
\label{rem:int_constraints_geo}
Just as a hierarchy
$\hier_1\in\nset\times\rset_+^{L+1}\times\rset_+^{L+1}$ solving
Problem~\ref{prob:cont_opt} must be adjusted to satisfy the practical
constraints of the discretization,
$\hier_1\approx\hier\in\nset\times\discreteh^{L+1}\times\nset_+^{L+1}$,
so must a hierarchy that is geometric with a general $\levsep$.
Hence, the restriction to general geometric sequences of mesh sizes,
without the true constraint
$\{h_0\levsep^{-\ell}\}_{\ell=0}^L\in\discreteh^{L+1}$, offers no
practical improvement over the more general optimization in Section~\ref{s:nongeo};
we merely include the comparison here to point out that one can often
find geometric hierarchies that are close to optimal hierarchies.
Figure~\ref{fig:int_work} shows the effect of applying these domain
constraints to the number of elements and number of samples on the
optimality of the hierarchies.  This figure compares the
work~\eqref{eq:model_work} of five hierarchies:
\begin{enumerate}
\item The ``real-valued'' optimized hierarchy with $h_\ell$ defined by
  \eqref{eq:hl_optimal} and $M_\ell$ defined by \eqref{eq:opt_M_ell},
\item The ``integer-valued'' hierarchy obtained by ceiling $M_\ell$ in
  \eqref{eq:opt_M_ell} and $h_\ell^{-1}$ in \eqref{eq:hl_optimal} to
  obtain an integer number of samples and an integer number of
  elements, respectively,
\item The real-valued geometric
  hierarchy with $\levsep$ as defined by \eqref{eq:opt_beta},
  $h_0=0.5$ and $M_\ell$ again as defined by
  \eqref{eq:opt_M_ell},
\item The integer-valued geometric hierarchy obtained by using the
  ceiling of both $M_\ell$ in \eqref{eq:opt_M_ell} and the previous
  $h_\ell^{-1}$.
\item Finally, a hierarchy obtained by performing a limited
  brute-force integer optimization in the neighboring integer space
  around the optimized $h_\ell^{-1}$ and $M_\ell$ in this work.
\end{enumerate} In all cases, we use the
parameters of \textbf{Ex.1} in Table~\ref{tbl:problem_params}.  Similar plots can be
produced with different values.  On the other hand, the number of
levels, $L$, was numerically optimized and chosen according to
Figure~\ref{fig:int_L}.  These plots show that simply taking the
ceiling of the number of samples and number of elements produces a
hierarchy that is nearly optimal.  Notice also in
Figure~\ref{fig:int_L} that the optimal $L$ of the optimized
real-valued hierarchies is well within the developed bounds
\eqref{eq:L_bound_TOL}, up to an integer rounding.  However, the
bounds no longer hold when considering integer-valued hierarchies.
\end{remark}

\begin{remark}[When optimal hierarchies become geometric]
Recall that for $\chi\neq1$, given $L$, $h_0$, and $h_L$ the optimal
intermediate mesh sizes satisfy~\eqref{eq:h_inner} which corresponds
to the ratios
\begin{align}
\label{eq:hl_ratio}
  \frac{h_\ell}{h_{\ell+1}} & =
  \left(\frac{h_0}{h_L}\right)^{\chi^\ell \frac{1-\chi}{1-\chi^L}}
  \chi^{\frac{2}{d\gamma}\left(  \frac{L \chi^{\ell}}{1-\chi^L} - \frac{1}{1-\chi}\right)}
\end{align}
between successive levels.
Keep the optimal choices of $h_0(\theta,L)$ and $h_L(\theta,L)$
in~(2.19a), let $\levsep$ satisfy~\eqref{eq:opt_beta}
and pick $\theta$ to solve the equation
\begin{align}
  \label{eq:force_geometric}
  h_L(\theta,L) & = h_0(\theta,L)\levsep^{-L}
\end{align}
given $L$. This means choosing the error splitting parameter
\[ \theta = 1-Q_W\tol^{-1}\levsep^{-q_1(L+1)}
\left(\frac{V_0}{Q_S}\right)^{\frac{\eta}{x}}\]
for any  $L \in \nset$ sufficiently large to make
the bias smaller than the total error tolerance, that is $\theta>0$.
Now, by~\eqref{eq:force_geometric} the ratios~\eqref{eq:hl_ratio}
reduces to $h_{\ell}/h_{\ell+1} = \beta$.
In other words, given $\tol$, to any sufficiently large $L$ there
corresponds at least one (suboptimal) $\theta$ such that the
corresponding optimized hierarchy is geometric.
\end{remark}

\section{Numerical results}\label{s:res}\makeatletter{}In this section, we first introduce the two test problems: a geometric
Brownian motion SDE for which $\chi>1$ and a random PDE for which
$\chi<1$ or $\chi>1$, depending on the linear solver used to solve the
corresponding linear system.  We then describe several implementation
details and finally conclude by presenting the actual numerical
results.  We do not show results for the case $\chi=1$ since we proved
that geometric hierarchies are optimal in this case and similar
results can be found in the standard work of Giles \cite{giles08}.

\subsection{Overview of examples} \label{sec:prob-overview}
We consider two numerical examples for which we can compute a reference solution.
\subsubsection{Ex.1}
This problem is based on Example~\ref{ex:spde_problem} in
Section~\ref{s:hier_intro} with some particular choices that satisfy
the assumptions therein.  First, we choose ${\mathcal{D} = (0,1)^3}$
and assume that the forcing is
\[ f(\mathbf x) = f_0 + \widehat f \sum_{i=0}^K \sum_{j=0}^K
\sum_{k=0}^K \Phi_{ijk}(\mathbf x) Z_{ijk},\] where
\[ \Phi_{ijk}(\mathbf x) = \sqrt{\lambda_i \lambda_j \lambda_k}
\phi_i(x_1)\phi_j(x_2)\phi_k(x_3), \] and
\begin{align*}
    \phi_i(x) &=
\begin{cases}
    \cos\left( \frac{5 \Lambda i}{2} \pi x \right), & i \text{ is even}, \\
    \sin\left( \frac{5 \Lambda (i+1)}{2} \pi x \right), & i \text{ is odd},
\end{cases} \\
\lambda_i &=  \left( 2 \pi \right)^\frac{7}{6}\Lambda^\frac{11}{6}
\begin{cases}
    \frac{1}{2},  & i=0, \\
    \exp\left( - 2\left( \pi \frac{i}{2} \Lambda \right)^2 \right), & i \text{ is even}, \\
    \exp\left( - 2\left( \pi \frac{i+1}{2} \Lambda \right)^2 \right), & i \text{ is odd},
\end{cases}
\end{align*}
for given parameters $\Lambda$, positive, and $K$, positive integer, and $\mathbf Z = \lbrace Z_{ijk} \rbrace$, a set of $(K+1)^3$ i.i.d. standard normal random variables.
Moreover, we choose the diffusion coefficient to be a function of two random variables as follows:
\begin{align}
    a(\mathbf x) &= a_0 + \exp \Big(4 Y_1 \Phi_{121}(\mathbf x) + 40 Y_2 \Phi_{877}(\mathbf x)\Big).
\end{align}
Here, $\mathbf Y = \lbrace Y_1, Y_2 \rbrace$ is a set of i.i.d. standard normal random variables, also independent of $\mathbf Z$.
Finally, we choose the quantity of interest, $g$, as a localized average around a point $\mathbf x_0$,
\[g = \left( 2 \pi \sigma^2 \right)^\frac{-3}{2} \int_\mathcal{D} \exp \left( - \frac{ \| \mathbf x - \mathbf x_0 \|^2_2}{2 \sigma^2} \right) u(\mathbf x) d\mathbf x, \]
and select the parameters $a_0=0.01, f_0=50, \widehat f=10, \Lambda = \frac{0.2}{\sqrt{2}}, K=10, \sigma^2=0.02622863$ and
${\mathbf x_0 = \left[0.5026695,0.26042876,0.62141498\right]}$.
Since the diffusion coefficient, $a$, is independent of the forcing, $f$, a reference solution can be calculated to sufficient accuracy
by scaling and taking expectation of the weak form with respect to $\mathbf Z$ to
obtain a formula with constant forcing for the conditional expectation
with respect to $\mathbf Y$.
We then use stochastic collocation, \cite{bnt2010}, with a sufficiently accurate quadrature to produce the reference value, $\E{g}$.
Using this method, the reference value $1.6026$ is computed with an error estimate of $10^{-4}$.
\subsubsection{Ex.2}
    The second example is a one-dimensional geometric Brownian motion based on Example~\ref{ex:sde_problem} where we make the following choices:
    \begin{align*}
        T &=1, \\
        a(t,u) &= 0.05 u,\\
        b(t,u) &= 0.2 u, \\
        g(u) &= 10 \max(u(1)-1,0).
    \end{align*}
    The exact solution can be computed using a change of variables and It\^o's formula.
    For the selected parameters, the solution is $1.04505835721856$.

\subsection{Implementation and runs}
To test the different hierarchies presented in this work we extend the
CMLMC algorithm~\cite{haji_CMLMC} to optimal hierarchies and implement
it in the C programming language. The CMLMC algorithm solves problems
with tolerances larger than the requested $\tol$ to cheaply get
increasingly accurate estimates of the constants $Q_W$ and $Q_S$ and
the variances $V_\ell$ for all $\ell = 0, 1 \ldots L$. This is
achieved in a Bayesian setting to incorporate the generated samples
with the models \eqref{eq:model_error}.

We stress that with these numerical results we aim to illustrate what
happens when the hierarchies of Section~\ref{s:opt} are used in a
practically viable algorithm which approximates several relevant
parameters during the computations. This means that we do not directly
observe the optimal hierarchies as derived in the continuous
optimization.

Unless $\theta$ is specified explicitly, our CMLMC algorithm uses a
\textit{computational} splitting parameter, $\theta$, calculated based
on the expected bias as
\begin{equation}
\label{eq:comp_theta}
\theta = 1 - \frac{Q_W h_L^{q_1}}{\tol},
\end{equation}
to relax the statistical error constraint. If $h_L$ satisfies
 \eqref{eq:hl_optimal_chi1} or \eqref{eq:hl_optimal}, this is the same as
the optimal splitting parameter defined by~\eqref{eq:opt_beta_chi1}
 or \eqref{eq:theta_opt}, respectively.

For implementing the solver for the PDEs in test problem {\bf Ex.1},
we use PetIGA~\cite{Collier2013}.  While the primary intent of
this framework is to provide high-performance B-spline-based finite
element discretizations, it is also useful in applications where the
domain is topologically square and subject to uniform refinements. As
its name suggests, PetIGA is designed to tightly couple to
PETSc~\cite{petsc-efficient}. The
framework can be thought of as an extension of the PETSc library,
which provides methods for assembling matrices and vectors that result
from integral equations.  We use uniform meshes with a standard
trilinear basis to discretize the weak form of the model problem,
integrating it with eight quadrature points. We also generate results
for two linear solvers for which PETSc provides an interface.  The
first solver is an {\bf Iterative} GMRES solver that solves a linear
system in almost linear time with respect to the number of degrees of
freedom for the mesh sizes of interest; in other words, in this case
$\gamma=1$ and $\chi>1$.  The second solver is the {\bf Direct} solver
MUMPS~\cite{Amestoy2001}. For the mesh sizes of interest,
the running time of MUMPS varies from quadratic to linear in the total
number of degrees of freedom. The best fit turns out to be
$\gamma=1.5$ in this case, which gives $\chi < 1$.  From
Corollary~\ref{thm:opt_work_chi_n1} (or
Corollary~\ref{thm:opt_geo_work}), the complexity for all the examples
is expected to be ${\Order{\tol^{-s}}}$, where $s$ depends on
$q_1, q_2,$ and $d\gamma$.  These and other problem parameters are
summarized in Table~\ref{tbl:problem_params} for the different
examples. Also included in this table is the optimal level separation
constant $\levsep$, which we used when computing with geometric
hierarchies.  Obviously, as mentioned in
Remark~\ref{rem:int_constraints_opt}, the ``real-valued'' hierarchies
we derived cannot always be used in practice and we follow the
strategies outlined in that remark to produce ``integer-valued''
hierarchies that can be used.

\begin{table}
\centering
\begin{tabular}{|r|c|c|c|}
  \hline
&  {\bf Ex.1} &  {\bf Ex.1} & {\bf Ex.2} \\ \hline
$d$ &  \multicolumn{2}{|c|}{3}  & 1\\  \hline
$q_1$ & \multicolumn{2}{|c|}{2} & 1\\ \hline
$q_2$ & \multicolumn{2}{|c|}{4} & 2\\ \hline
Estimated $V_0$ & \multicolumn{2}{|c|}{0.0565} & 1.7805 \\ \hline
Estimated $Q_W$ & \multicolumn{2}{|c|}{1.3653} & 0.0307 \\ \hline
Estimated $Q_S$ &\multicolumn{2}{|c|}{0.1519} & 0.2630\\ \hline
Solver & GMRES & MUMPS & Milstein \\ \hline
$\gamma$ & 1 & 1.5 & 1\\  \hline
$\chi$ & 4/3 & 8/9 & 2 \\ \hline
$\eta$ & 2/3 & 4/9 & 1 \\ \hline
$s$ & 2 & 2.25& 2\\ \hline
Optimal $\levsep$ & 1.7778 & 1.6018 & 4\\ \hline
Work to seconds constant & $10^{-4}$ &  $3\times 10^{-6}$ &  $9 \times
                                                            10^{-8}$ \\\hline
\end{tabular}
\caption{Summary of problem parameters.}
\label{tbl:problem_params}
\end{table}

We run each setting $100$ times and show in plots in the next section the medians with vertical bars spanning from the $5\%$ percentile
to the $95\%$ percentile.
Finally, all results were generated on the same machine with $52$ gigabytes of memory to ensure that no overhead is introduced due to hard disk access
during swapping that could occur when solving the three-dimensional PDEs with a fine mesh.
We use the parameters listed in Table~\ref{tbl:alg_params} for the CMLMC algorithm \cite{haji_CMLMC}.

\begin{table}
\centering
\begin{tabular}{|p{0.75in}|p{2in}|p{0.9in}|p{0.9in}|}
\hline
Parameter & Purpose & Value for {\bf Ex.1} & Value for {\bf Ex.2} \\
\hline
$\kappa_0$ and $\kappa_1$ & Confidence parameter for the weak and strong error models &  $0.1$ for both & $0.1$ for both\\
\hline
$\tol_{\max}$ &The maximum tolerance with which to start the algorithm. & 0.5 & 0.1\\
\hline
$r_1$ and $r_2$ &Controls computational burden to calibrate the problem parameters. & $2$ and $1.1$, respectively & $2$ and $1.1$, respectively \\
\hline
Initial hierarchy &The initial hierarchy to start the CMLMC algorithm.& $L=2$ and $h_\ell^{-1} = \{4,6,8\}$ and $M_\ell = 10$ for all $\ell$. & $L=2$ and $h_\ell^{-1} = \{1,2,4\}$ and $M_\ell = 10$ for all $\ell$.  \\
\hline
$L_{\text{inc}}$ &Maximum number of values to consider when optimizing for L. & 2 & 2 \\
\hline
$\mathfrak{L}$ &Maximum number of levels used to compute estimates of $Q_W$ and $Q_S$. & 3 & 5  \\
\hline
$C_\alpha$ & Parameter related to the confidence in the statistical constraint &2 & 2\\
\hline
\end{tabular}
\caption{Summary of parameter values used in the CMLMC algorithm in our numerical tests. This table is reproduced from \cite{haji_CMLMC} where more information is available.}
\label{tbl:alg_params}
\end{table}

\subsection{Results}
\begin{figure}
\centering
\includegraphics[scale=0.5]{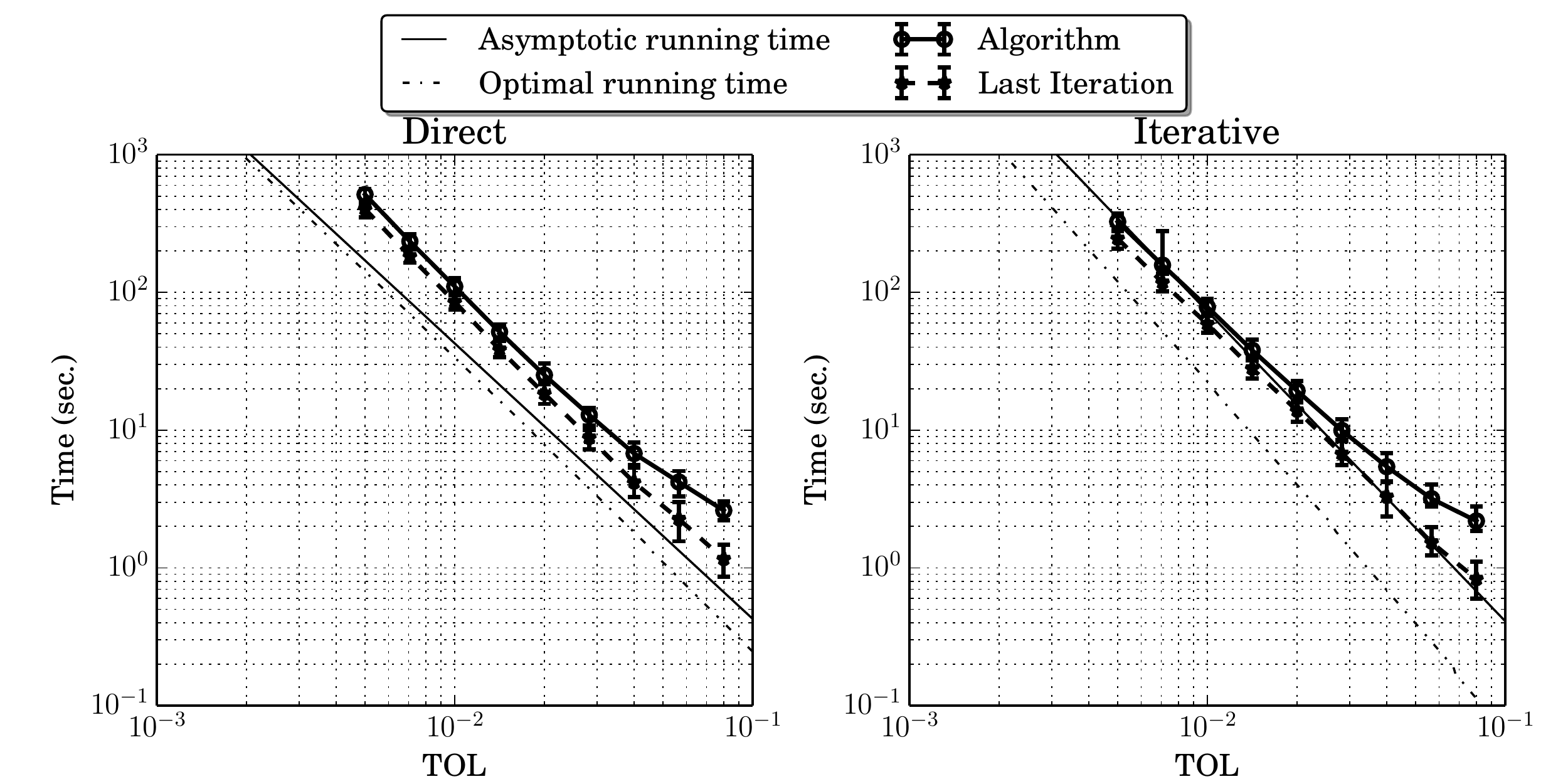}
\caption{\textbf{Ex.1}: The running time of the CMLMC algorithm and
  its last iteration when using optimal hierarchies.  Using the
  parameters in \ref{tbl:problem_params}, we also show in this figure
  the expected running time when using the optimal, unconstrained
  hierarchy defined by \eqref{eq:parameters} and the expected
  asymptotic work according to Corollary~\eqref{thm:opt_work_chi_n1}.}
\label{fig:ex1-runtime-rate}
\end{figure}

\begin{figure}
 \centering
 \includegraphics[scale=0.5]{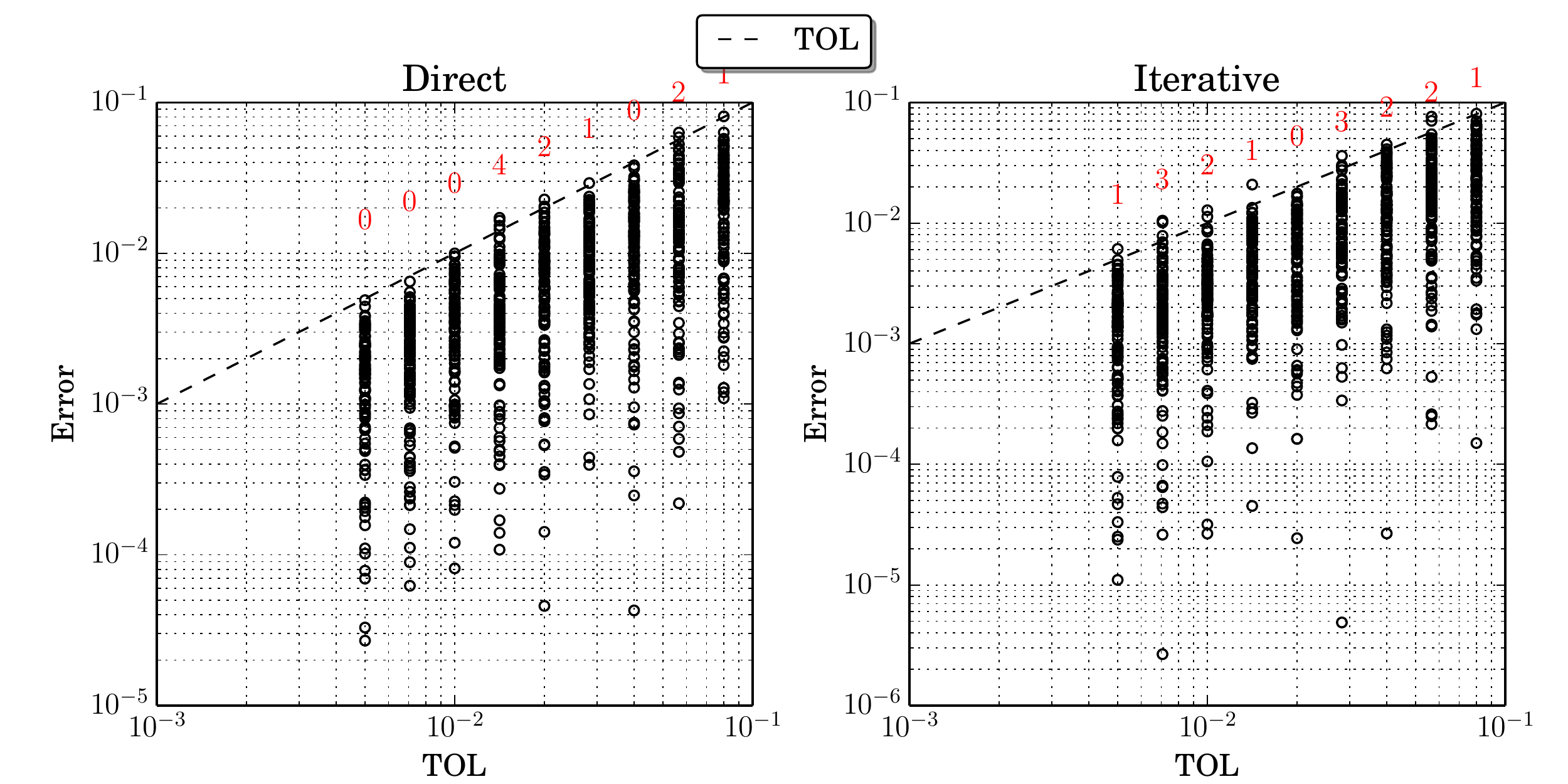}
 \includegraphics[scale=0.5]{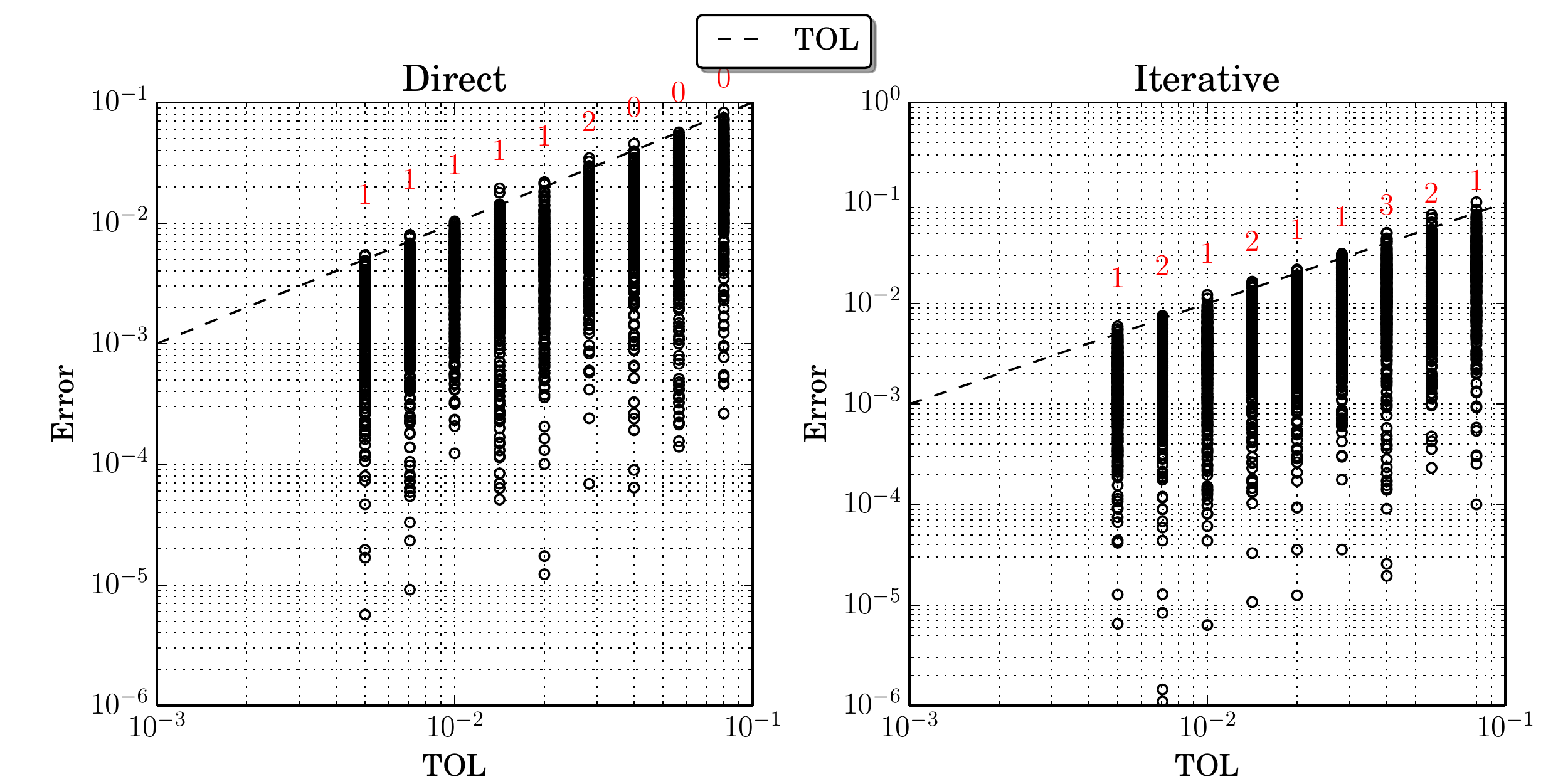}
 \caption{\textbf{Ex.1}: The true errors calculated using optimal
   hierarchies in the top plot, and using geometric hierarchies in the
   bottom one. Each point shows the true error of an independent run
   of MLMC.  The numbers on top of the $\tol$ line are the percentage
   of algorithm runs that produced a larger error than the required
   tolerance. Notice that the choice $C_\alpha=2$ gives a confidence
   of a $95\%$ in the error bound, as predicted in \cite[Lemma
   A.2]{haji_CMLMC}. }
 \label{fig:ex1-ex-errors}
 \end{figure}

\begin{figure}
 \centering
 \includegraphics[scale=0.5]{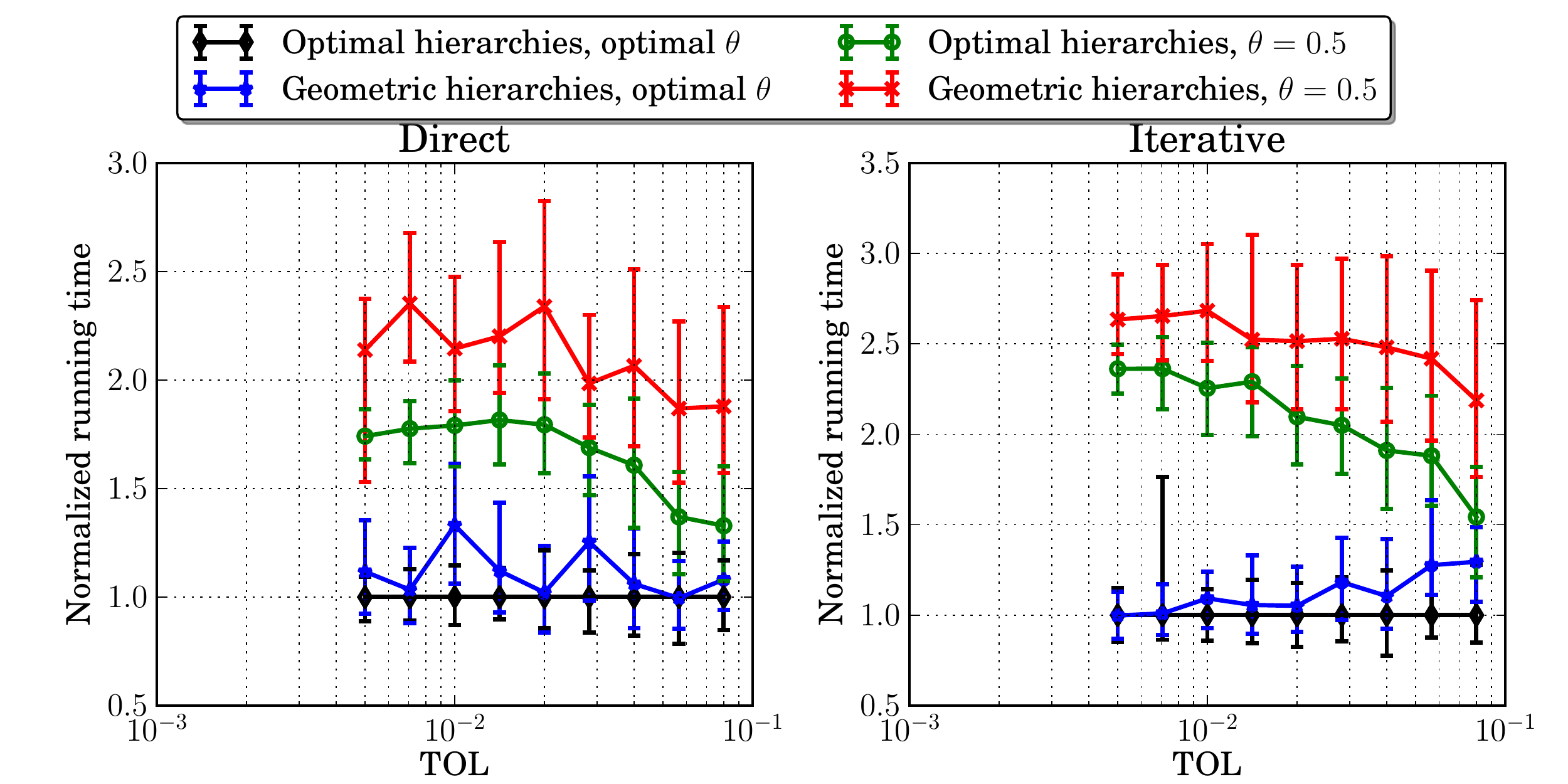}
 \caption{\textbf{Ex.1}: Actual running time of the CMLMC algorithm
   when using optimal and geometric hierarchies with different
   tolerance splitting, normalized by the average running time of the
   algorithm when using optimal hierarchies. Compare this figure to
   Figure~\ref{fig:int_work}, where the latter is based on the
   theoretical results.  Observe that most of the gain in
   computational complexity is due to the choice~\eqref{eq:comp_theta}
   of $\theta$
   and that using optimal hierarchies does not significantly improve the
   running time over geometric hierarchies.}
 \label{fig:ex1-runtime-const}
 \end{figure}

\begin{figure}
 \centering
 \includegraphics[scale=0.5]{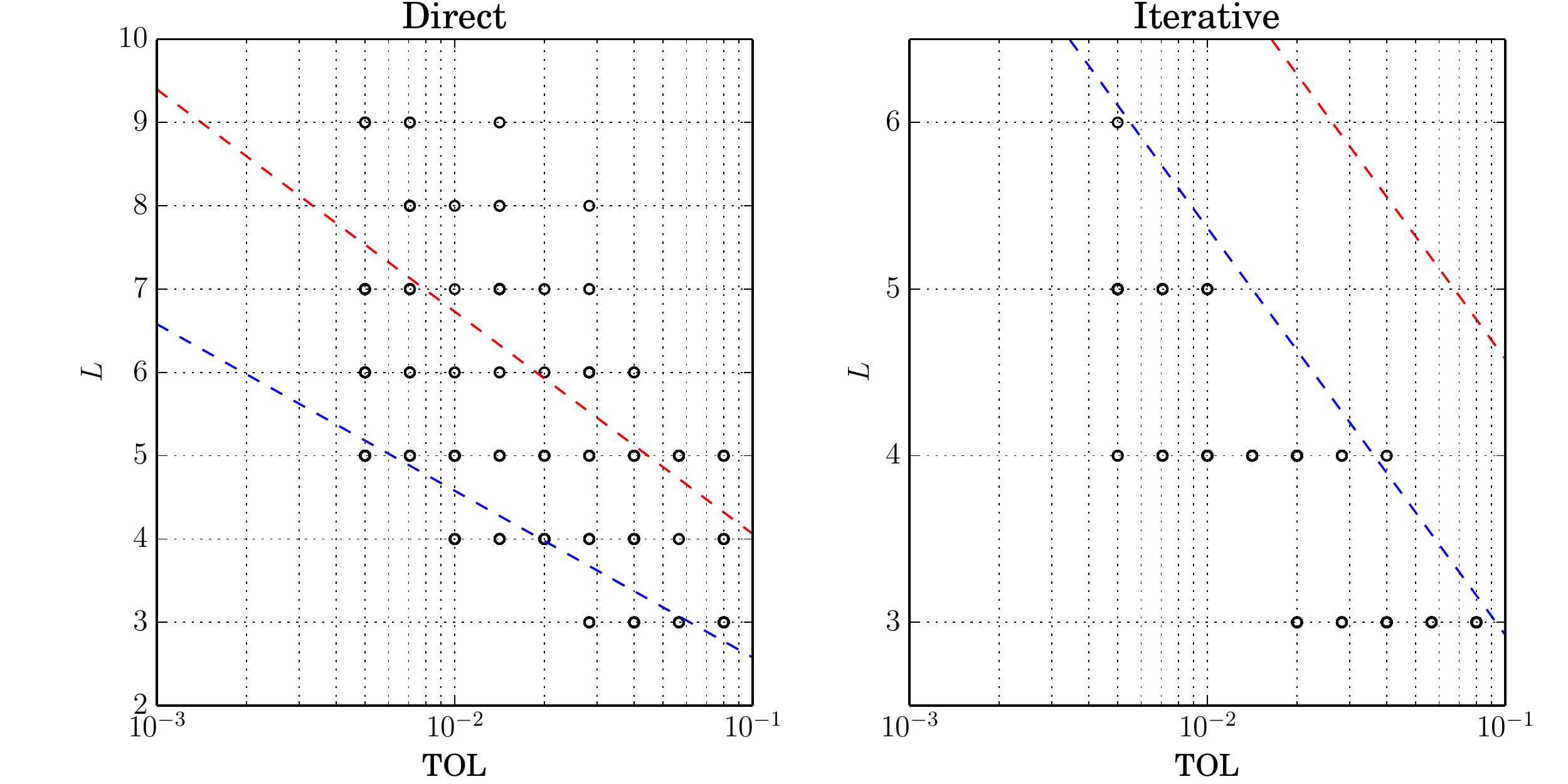}
 \caption{\textbf{Ex.1}: The used number of levels, $L$, for different
   tolerances in the last iteration of the CMLMC algorithm when using
   optimal hierarchies and ceiling $h_\ell^{-1}$ and $M_\ell$.
   Here, different circles correspond to independent runs of CMLMC.
   Compare this figure to the Figure~\ref{fig:int_L}, where the latter
   is based on the theoretical results.  The bounds are taken from
   \eqref{eq:L_bound_TOL}. The $L$ values used by the CMLMC algorithm
   fall outside the predicted bounds because the bounds are valid for
   the real-valued optimal hierarchies only.  On the other hand, CMLMC
   is constrained to the discrete sets of feasible hierarchies and further
   limits the increments of $L$ across iterations.  }
 \label{fig:ex1-L-vs-tol}
 \end{figure}

\begin{figure}
 \centering
 \includegraphics[scale=0.5]{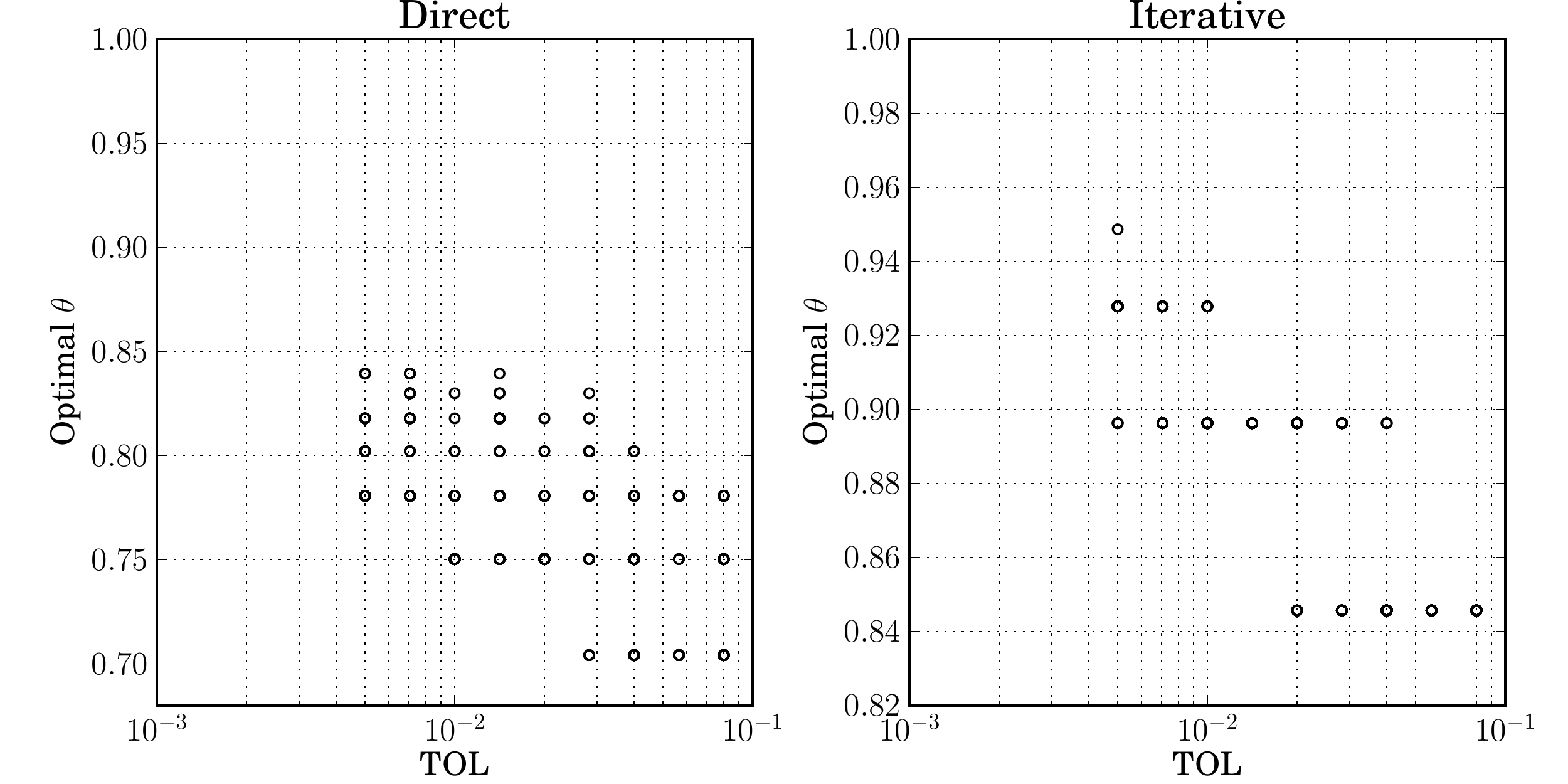}
 \caption{\textbf{Ex.1}: The error splitting parameter, $\theta$, as
   defined by \eqref{eq:theta_opt}. Recall that $\theta$ depends on
   $L$ and notice the correspondence to $L$ values in Figure~\ref{fig:ex1-L-vs-tol}}
 \label{fig:ex1-theta-vs-tol-exp}
 \end{figure}

\begin{figure}
 \centering
 \includegraphics[scale=0.5]{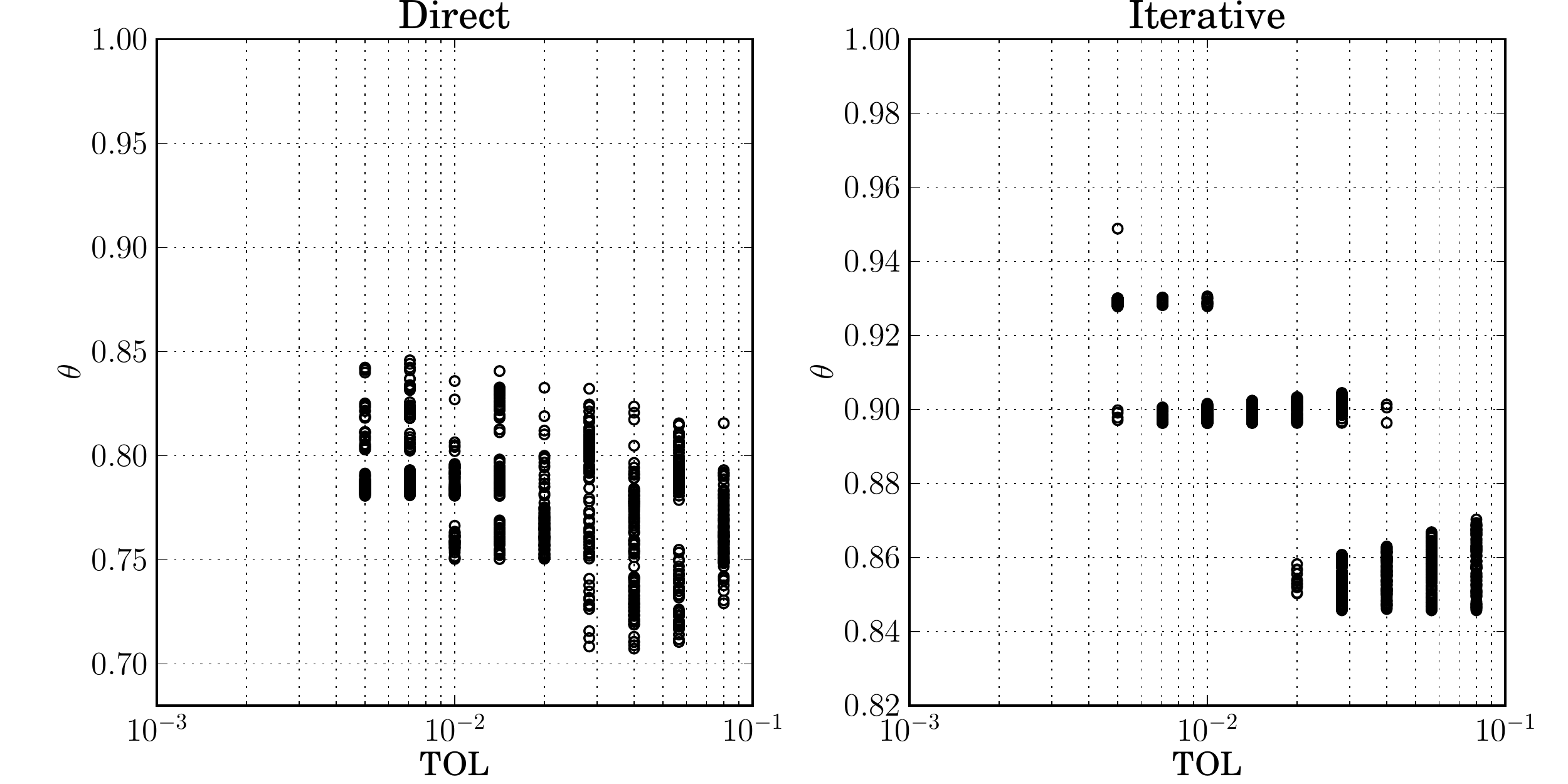}
 \caption{\textbf{Ex.1}: The computational splitting parameter,
   $\theta$, as defined in \eqref{eq:comp_theta} for the CMLMC
   algorithm for optimal hierarchies. Though these $\theta$ values
   correspond to the value in Figure~\ref{fig:ex1-theta-vs-tol-exp}, the
   differences are due to enforcing constraints on $h_L$.}
 \label{fig:ex1-theta-vs-tol-opt}
 \end{figure}

\begin{figure}
 \centering
 \includegraphics[scale=0.5]{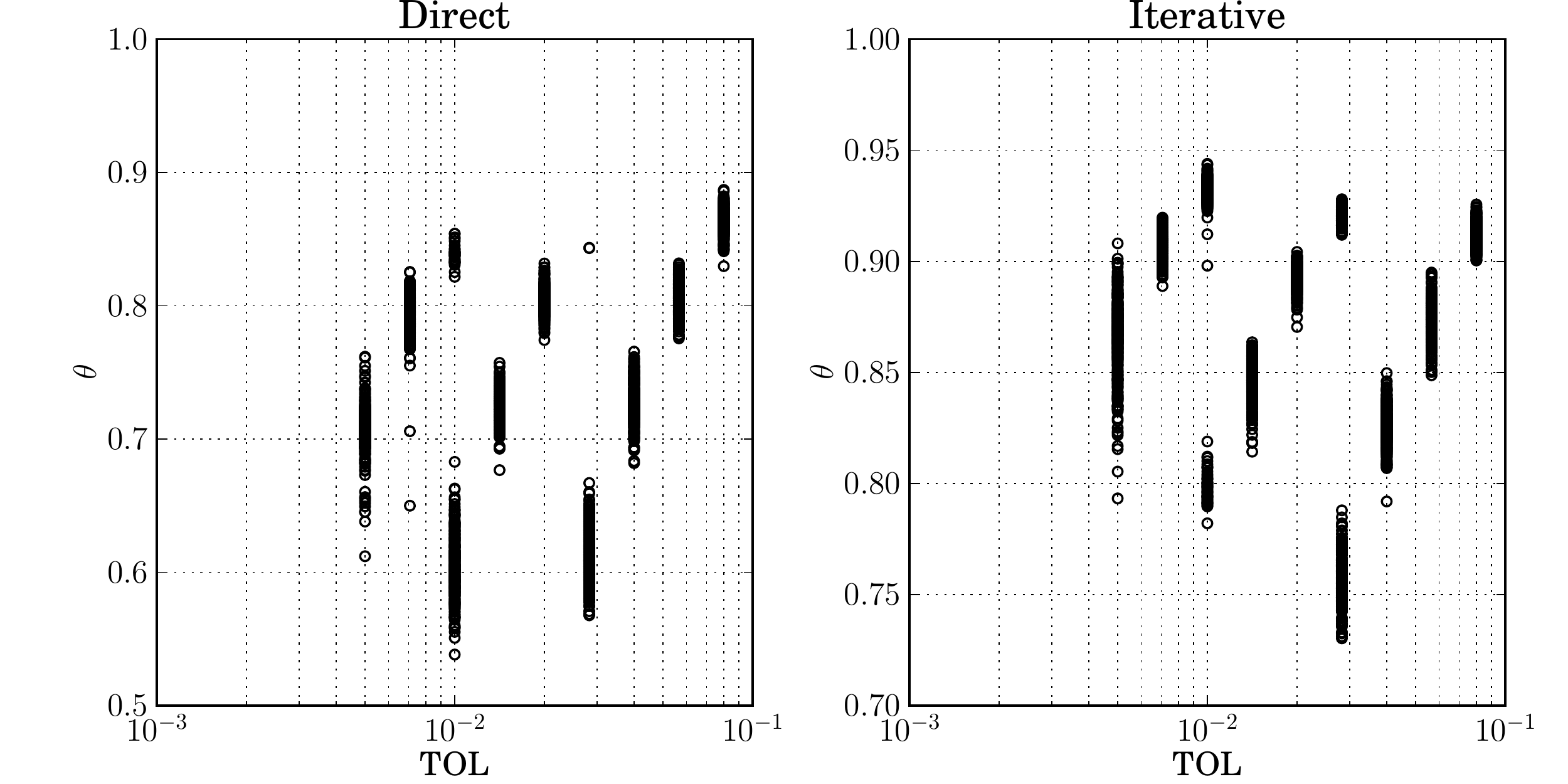}
 \caption{\textbf{Ex.1}: The computational splitting parameter, $\theta$, in the CMLMC algorithm for optimal geometric hierarchies.}
 \label{fig:ex1-theta-vs-tol-geo}
 \end{figure}

\begin{figure}
 \centering
 \includegraphics[scale=0.5]{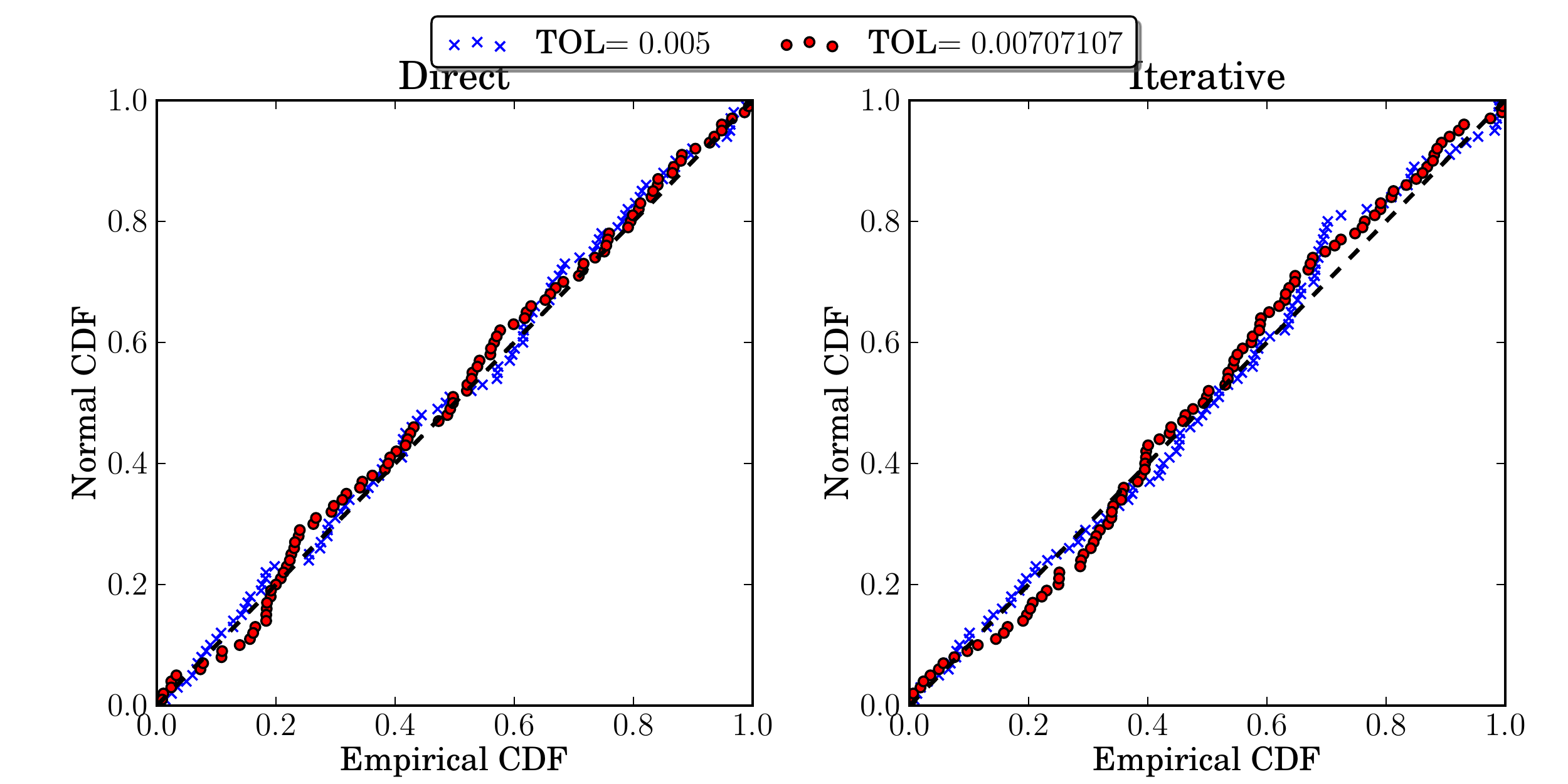}
 \caption{\textbf{Ex.1}: A QQ-plot indicating that, even when using non-geometric hierarchies, the distribution of the normalized statistical error is well approximated by the standard normal density.
  The work \cite[Lemma A.2]{haji_CMLMC} proved such results for geometric hierarchies.}
 \label{fig:ex1-qq-plot}
 \end{figure}

We start by presenting the results of \textbf{Ex.1}.  We show in
Figure~\ref{fig:ex1-runtime-rate} the total running time of the CMLMC
algorithm and its last iteration when using optimal hierarchies, after
taking the ceiling of the optimal number of elements $h_\ell^{-1}$ in
\eqref{eq:hl_optimal} and the optimal number of samples $M_\ell$ in
\eqref{eq:opt_M_ell}.  Using the parameters in
\ref{tbl:problem_params}, we also show in this figure the expected
running time when using the optimal, unconstrained hierarchy defined
by \eqref{eq:parameters} and the expected asymptotic work according to
Corollary~\eqref{thm:opt_work_chi_n1}.  This figure shows good
agreement between the expected theoretical results and the actual
final running time.  Figure~\ref{fig:ex1-ex-errors} shows that the
true error that was computed using the reference solution when using
optimal hierarchies is less than the required tolerance with the
required confidence of $95\%$, in accordance with the chosen value of
$C_\alpha=2$ and \cite[Lemma A.2]{haji_CMLMC}.
Figure~\ref{fig:ex1-runtime-const} compares the computational
complexity of optimal hierarchies to geometric hierarchies for
different values of $\theta$.  This figure shows numerical
confirmation that optimal hierarchies do not give significant
improvement over geometric hierarchies, especially for optimal values
of $\theta$.  In other words, the improvement of the running time is
mainly due the better choice of $\theta$ as discussed in
\cite{haji_CMLMC}.  Figure~\ref{fig:ex1-theta-vs-tol-exp} shows the
optimal splitting, $\theta$, as defined by \eqref{eq:theta_opt}.
Compare this figure to Figure~\ref{fig:ex1-L-vs-tol}, which shows the
used number of levels, $L$, for different tolerances, and notice the
dependence of $\theta$ on the number of levels, $L$.  On the other
hand, Figure~\ref{fig:ex1-theta-vs-tol-opt} shows the computational
splitting used in the CMLMC algorithm.  Notice that $\theta$ follows a
similar pattern in both Figure~\ref{fig:ex1-theta-vs-tol-exp} and
Figure~\ref{fig:ex1-theta-vs-tol-opt}. The continuous change in the
latter is due to differences in the estimation of $Q_W$ for different
runs of the algorithm.  For comparison,
Figure~\ref{fig:ex1-theta-vs-tol-geo} shows that the computational
splitting parameter produced when using geometric hierarchies is
different from the computational splitting parameter produced when
using optimal hierarchies. However, as $\tol \to 0$, they both seem to
be not too far from the limit in \eqref{eq:theta_opt_A}.
Finally, even though \cite[Lemma A.2]{haji_CMLMC} assumes a geometric sequence, Figure~\ref{fig:ex1-qq-plot} shows that the lemma still holds for non-geometric hierarchies;
    i.e., that the cumulative density function (CDF) of the true error when suitably normalized is well approximated by a standard normal CDF.

\begin{figure}
\centering
\includegraphics[scale=0.5]{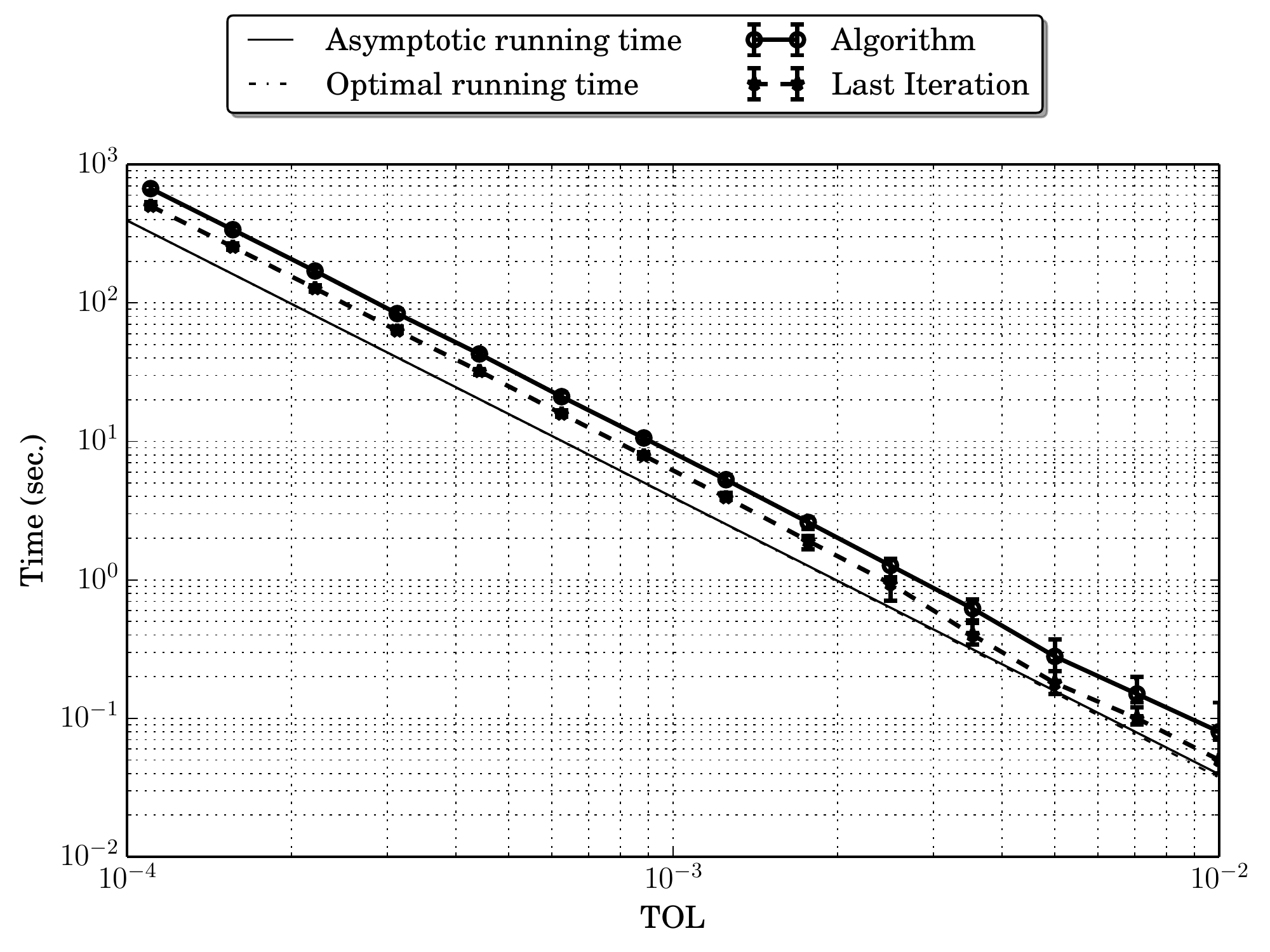}
\caption{\textbf{Ex.2}: The running time of the CMLMC algorithm. The reference dashed line is $\Order{\tol^{-2}}$ as predicted in \cite[Theorem 2.3]{tsgu13}.}
\label{fig:ex2-runtime-rate}
\end{figure}

\begin{figure}
 \centering
 \includegraphics[scale=0.5]{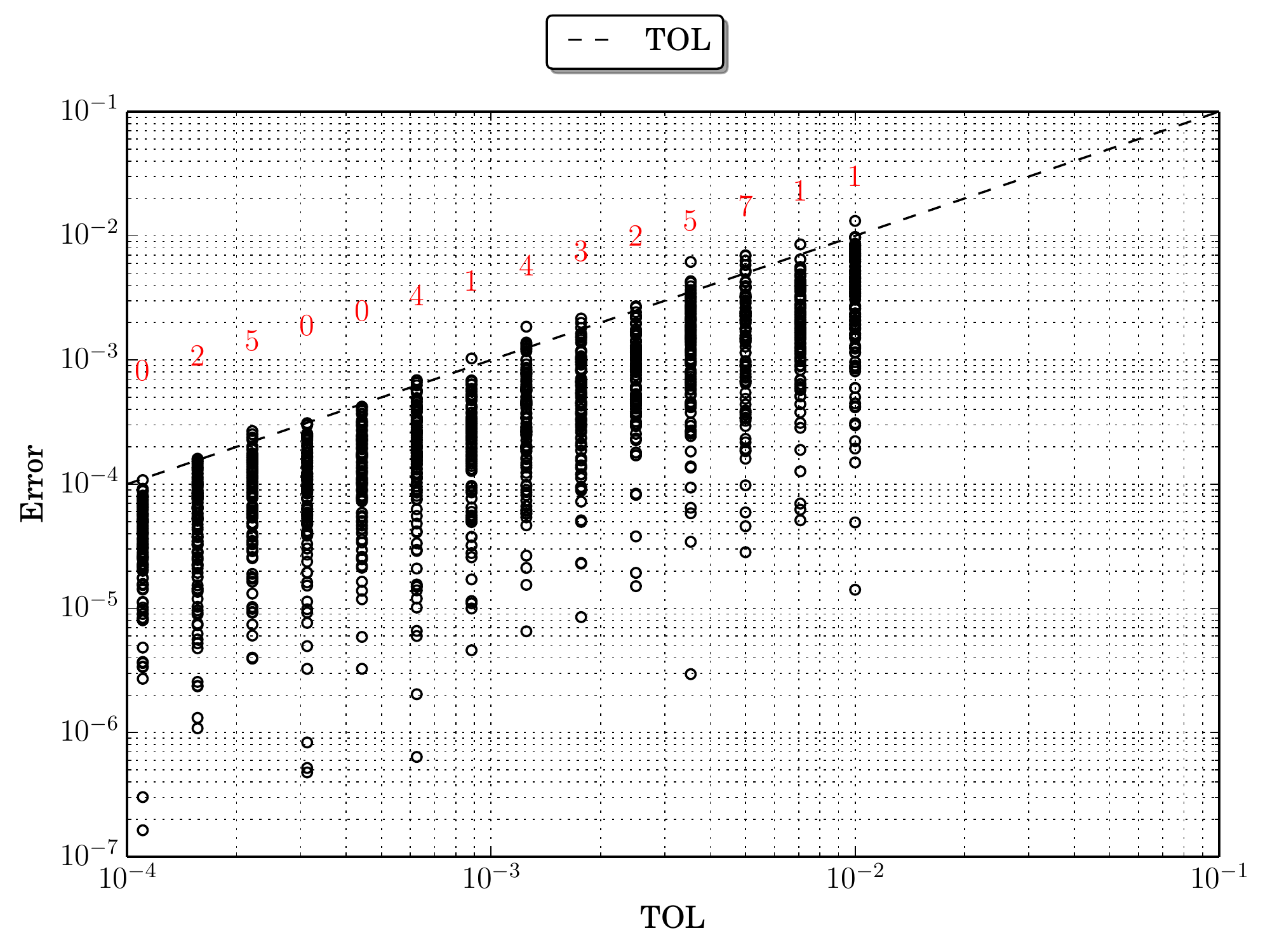}
 \caption{\textbf{Ex.2}: The true errors calculated using geometric hierarchies with $\levsep=4$.
     The numbers on top of the $\tol$ line are the percentage of algorithm runs that produced a larger error than the required
         tolerance. Remember that $C_\alpha=2$ gives a confidence of a $95\%$ in the error bound, as predicted in \cite[Lemma A.2]{haji_CMLMC}.}
 \label{fig:ex2-ex-errors}
 \end{figure}

Next, we focus on \textbf{Ex.2} where $\chi=2$ using the Milstein scheme. Since we showed previously that geometric hierarchies are near-optimal, we only present the results when using geometric hierarchies in this case.
The optimal geometric constant, $\levsep$, is $4$ in this case according to \eqref{eq:opt_beta}.
Figure~\ref{fig:ex2-runtime-rate} shows that the actual running time
of the CMLMC algorithm has the expected rate $\tol^{-2}$, again as
predicted in Corollary~\ref{thm:opt_geo_work}, or indeed for these
nested geometric hierarchies already in~\cite{giles08_m}.
Figure~\ref{fig:ex2-ex-errors} shows that the true errors for different tolerances are less than the required tolerance with the required confidence of $95\%$.

\section{Conclusions}\label{s:conc}\makeatletter{}

MLMC sampling methods are becoming increasingly popular due to their
robustness and simplicity.  In this work, in
Theorems~\ref{thm:opt_hier_chi_1} and~\ref{thm:opt_hier_chi_n1} and
Corollary~\ref{thm:opt_geo_work}, we have developed optimal
non-geometric and geometric hierarchies for MLMC by assuming certain
asymptotic models on the weak and strong convergence and the average
computational cost per sample.  While it is important to optimize the
geometric level separation parameter, $\beta$, and the tolerance
splitting parameter, $\theta$, to obtain significant computational
savings, we have shown, in Remark~\ref{rem:asymb_work}, that with
these optimized parameters, geometric hierarchies are
nearly optimal and that, asymptotically, their computational
complexity is the same as the non-geometric optimal hierarchies.
Moreover, we have analyzed the asymptotic behavior of the optimal
tolerance splitting parameter, $\theta$, between the bias and the
statistical error contribution.  We have also discussed how enforcing
constraints on parameters of MLMC hierarchies affects the optimality
of these hierarchies.  These constraints include an upper and lower
bound on the mesh size or enforcing that the number of samples and the
number of discretization elements are integers.  Our numerical results
show remarkable agreement between our theory of optimal hierarchies
and their asymptotic behavior and the performance of the CMLMC
algorithm.

                             In future work, it is possible to improve the efficiency of the
     MLMC method by including certain non-asymptotic terms in the
     models for the weak and strong convergence or the computational
     complexity.  Moreover, since the asymptotic dependence of the
     computational complexity on the different problem constants is
     clearly shown in Corollaries~\ref{thm:opt_work_chi_1}
     and~\ref{thm:opt_work_chi_n1}, one can devise methods to combine
     with MLMC to reduce the total computational complexity by
     affecting these constants, for example by reducing the variance,
     $V_0$, for the case where $\chi>1$.

\section*{Acknowledgments}
R.~Tempone is a member of the Research Center on Uncertainty
Quantification (SRI-UQ), division of Computer, Electrical and
Mathematical Sciences and Engineering (CEMSE) at King Abdullah
University of Science and Technology (KAUST).  The authors would like
to recognize the support of the following KAUST and University of
Texas at Austin AEA projects: Round 2, ``Predictability and
Uncertainty Quantification for Models of Porous Media'', and Round 3,
``Uncertainty quantification for predictive modeling of the
dissolution of porous and fractured media''.  F. Nobile has been
partially supported by the Swiss National Science Foundation under the
Project No. 140574 ``Efficient numerical methods for flow and
transport phenomena in heterogeneous random porous media'' and by the
Center for ADvanced MOdeling Science (CADMOS). E. von~Schwerin has
been partially supported by the aforementioned SRI-UQ and CADMOS. We
would also like to acknowledge the following open source software
packages that made this work possible: {\tt
  PETSc}~\cite{petsc-web-page}, {\tt PetIGA}~\cite{Collier2013}.
\newpage
\appendix
\section{Derivations and proofs}\label{s:app}\makeatletter{}

\subsection{Optimal hierarchies given $h_0$, $\theta$, and $L$}
\label{s:app_opt_hier}
Here we solve Problem~\ref{prob:cont_opt} of Section~\ref{s:nongeo}
for the optimal hierarchy for any fixed value of $L$.
We initially treat the parameter $\theta$ as given, postponing
its optimization until later, and proceed in two steps to find the
optimal $\{M_\ell\}_{\ell=0}^L$ and $\{h_\ell\}_{\ell=0}^L$.
Assuming general work estimates $\{W_\ell\}_{\ell=0}^L$
in~\eqref{eq:work_sum} and general variance estimates of $\{V_\ell\}_{\ell=0}^L$,
 we assume equality in~\eqref{eq:var_requirement} and introduce the Lagrange
multiplier $\lambda$ to obtain the Lagrangian
\begin{align*}
    \mathcal{L}\left(\{M_\ell\}_{\ell=0}^L,\lambda\right) & =
  \sum_{\ell=0}^{L} M_\ell W_\ell
  +\lambda\left\{\sum_{\ell=0}^L \frac{V_\ell}{M_\ell}-
    \left(\theta\frac{\tol}{C_\alpha}\right)^2
  \right\}.
\end{align*}
The requirement that the variation of the Lagrangian with respect to
$M_\ell$ is zero, gives
  $M_\ell = \sqrt{\lambda\frac{V_\ell}{W_\ell}}.$
Solving for $\lambda$ in the variance constraint \eqref{eq:var_requirement}
and substituting back leads to \eqref{eq:opt_M_ell}.
Substituting this optimal $M_\ell$ in the total work \eqref{eq:work_sum} yields
\begin{align}
  \label{eq:total_work_app}
  \work(\hier) & =
  \left( \frac{C_\alpha}{ \theta\tol} \right)^2
  \left( \sum_{\ell=0}^L \sqrt{W_\ell V_\ell} \right)^2.
\end{align}
We proceed to find the optimal $\{h_\ell\}_{\ell=0}^L$ under the particular models~\eqref{eq:cont_opt}.
The total work \eqref{eq:total_work_app} is minimized when
\begin{align}
  \label{eq:tmp_target}
  \sum_{\ell=0}^L \sqrt{W_\ell V_\ell} & =
  \sqrt{\frac{\vg}{h_0^{d\gamma}}} +
  \sqrt{Q_S}\sum_{\ell=1}^L \sqrt{\frac{h_{\ell-1}^{q_2}}{h_\ell^{d\gamma}}},
\end{align}
is minimized.
Here the finest mesh, $h_L$, is given by the bias
constraint~\eqref{eq:bias_constraint} as
\begin{align}
  \label{eq:h_L}
  h_L & = \left( \frac{(1-\theta)\,\tol}{Q_W} \right)^{\frac{1}{q_1}},
\end{align}
independently of the multilevel construction.
Now, treat the coarsest mesh, $h_0$, as given and find the optimal
$h_1,\ldots,h_{L-1}$ that minimize
\begin{align}
  \label{eq:work_factor_inner}
  \frac{1}{\sqrt{Q_S}}\sum_{\ell=1}^L \sqrt{W_\ell V_\ell} & =
  \sum_{\ell=1}^L \sqrt{\frac{h_{\ell-1}^{q_2}}{h_\ell^{d\gamma}}}.
\end{align}
The requirement that the derivative of this sum with respect to
$h_\ell$ equals zero, for $\ell=1,\ldots,L-1$, leads to the optimality
condition
\begin{align*}
  q_2 h_\ell^{\left(\frac{q_2+d\gamma}{2}\right)} & =
  d\gamma h_{\ell-1}^{\left(\frac{q_2}{2}\right)}
  h_{\ell+1}^{\left(\frac{d\gamma}{2}\right)},
\end{align*}
which after taking the logarithm and using $\chi$ defined
in~\eqref{eq:def_chi_eta}, leads to
\begin{align}
  \label{eq:lin_diff_eqn}
  - \log{\left(h_{\ell+1}\right)}
  + \left(1+\chi\right) \log{\left(h_\ell\right)}
  - \chi \log{\left(h_{\ell-1}\right)}
  & = -\frac{2}{d\gamma} \log{\left(\chi\right)}.
\end{align}
This is a second order linear difference whose solution depends on $\chi$.

\subsubsection{For $\chi=1$}
\label{s:app_chi_1}

This section provides proofs of Theorem~\ref{thm:opt_hier_chi_1}, Lemma~\ref{thm:opt_L_chi_1}, and Corollary~\ref{thm:opt_work_chi_1}.
The solution of the difference equation \eqref{eq:lin_diff_eqn} for the case $\chi=1$
 is the geometric sequence
\begin{align}
  \label{eq:h_inner_chi_1}
  h_\ell & = h_0\levsep^{-\ell}, && \text{with $\levsep=\left(\frac{h_0}{h_L}\right)^{1/L}$.}
\end{align}

In other words, all $h_\ell$ are defined in terms of $h_0$ and $h_L$, where the
latter is determined by $\theta$ through~\eqref{eq:h_L}
and we solve for the former by setting the derivative of~\eqref{eq:tmp_target}
with respect to $h_0$ equal to zero. This optimality condition becomes (for
        $q_2=d\gamma$)
\[ h_1 = \left(\frac{Q_S}{V_0}\right)^{\frac{1}{q_2}} h_0^2. \]
Combining this expression with \eqref{eq:h_inner_chi_1} for $\ell=1$ and solving for
$h_0$ yields
\begin{equation}
\label{eq:h0_chi1}
h_0 = h_L^\frac{1}{L+1} \left(\frac{V_0}{Q_S}\right)^\frac{L}{q_2 (L+1)}.
\end{equation}
Substituting this expression and \eqref{eq:h_L} in the expression for $\levsep$
in \eqref{eq:h_inner_chi_1} we obtain \eqref{eq:opt_beta_chi1}.
Moreover, substituting \eqref{eq:opt_beta_chi1} and \eqref{eq:h_inner_chi_1} and \eqref{eq:var_model} and \eqref{eq:wl_model} in
\eqref{eq:opt_M_ell} yields \eqref{eq:Ml_optimal_chi1}.
Next, we substitute \eqref{eq:Ml_optimal_chi1} and \eqref{eq:hl_optimal_chi1} in \eqref{eq:model_work} to obtain the optimal work for $q_2=d\gamma$
\begin{equation}
\label{eq:work_geo_chi1}
W = \left( \frac{C_\alpha}{\theta \tol} \right)^2 \left(\sqrt{V_0} h_0^\frac{-q_2}{2} + \sqrt{Q_S} \levsep^\frac{q_2}{2} L \right)^2.
\end{equation}
Using \eqref{eq:h_inner_chi_1} and \eqref{eq:h0_chi1}, we obtain
\begin{equation}
W = \left( \frac{C_\alpha}{\theta \tol} \right)^2 h_L^\frac{-q_2}{L+1} V_0^\frac{1}{L+1} Q_S^\frac{L}{L+1}\left(1+L\right)^2 .
\end{equation}
Substituting for $h_L$ from \eqref{eq:h_L}
\begin{equation}
W = \left(\frac{C_\alpha}{\theta \tol} \right)^2 \left( \frac{Q_W}{(1-\theta)\,\tol} \right)^{\frac{1}{ \eta(L+1)}} V_0^\frac{1}{L+1} Q_S^\frac{L}{L+1}\left(1+L\right)^2.
\end{equation}
Optimizing for $\theta$ yields \eqref{eq:opt_theta_chi1}. Substituting back gives the work as a function of $L$
\begin{equation}
W(L) = C_\alpha^2 \tol^{-2(1+\mathfrak{e}(L))}  Q_W^{2\mathfrak{e}(L)} V_0^{2 \eta \mathfrak{e}(L)} Q_S^{-2 \eta \mathfrak{e}(L)} Q_S \left(\frac{1}{2\eta}\right)^2 \left(1
        + \frac{1}{\mathfrak{e}(L)}\right)^{2(1+\mathfrak{e}(L))},
\end{equation}
where $\mathfrak{e}(L) = \frac{1}{2\eta(L+1)}$. Treating $L$ as a continuous variable and differentiating with respect to $L$ yields
\begin{equation}
\label{eq:WL_derv_chi1}
W'(L) = 2 W(L) \mathfrak{e}'(L) \left( C -y + \log(1+y) \right),
\end{equation}
where
${y=2\eta(L+1)} > 0$
and ${C = \log\left(\tol^{-1} Q_W V_0^\eta Q_S^{-\eta} \right)}$.
Setting \eqref{eq:WL_derv_chi1} to zero gives the equation
\begin{align}
  \label{eq:y_vs_C}
  y - \log(1+y) = C.
\end{align}
It follows that $\lim_{C \to \infty} \frac{y}{C} = 1$ which leads
to~\eqref{eq:L_chi_1_asymb} for the value of $L$
and~\eqref{eq:opt_asymp_work_chi1} for the work.
Since $y>0$, equation~\eqref{eq:y_vs_C} implies
\begin{equation}
  \label{eq:low_bound_y_over_C}
  1<y C^{-1}
\end{equation}
for any $C>0$.
Furthermore, for any $y>0$, it holds
\[ \left(\frac{\exp(1)-1}{\exp(1)}\right)(1+y) - C \leq 1+y -
\log(1+y) - C, \]
which together with~\eqref{eq:y_vs_C} gives
\begin{equation}
  \label{eq:upper_bound_y_over_C}
  y C^{-1} \leq \frac{\exp(1)}{\exp(1)-1}
    + \frac{1}{(\exp(1)-1)C},
\end{equation}
for any $C>0$. Inequalities~\eqref{eq:low_bound_y_over_C}
and~\eqref{eq:upper_bound_y_over_C} are~\eqref{eq:L_chi_1_bounds}.

\subsubsection{For $\chi \neq 1$}
\label{s:app_chi_n1}
This section provides proofs of Theorem~\ref{thm:opt_hier_chi_n1}, Lemma~\ref{thm:opt_L_chi_n1}, and Corollary~\ref{thm:opt_work_chi_n1}.
The solution of the difference equation \eqref{eq:lin_diff_eqn} for the case $\chi\neq 1$ is
\begin{align}
  \label{eq:h_inner}
  h_\ell & =
  h_0^{\left(\frac{\chi^\ell-\chi^L}{1-\chi^L}\right)}
  h_L^{\left(\frac{1-\chi^\ell}{1-\chi^L}\right)}
  \chi^{-\frac{2}{d\gamma}\left(\frac{L(1-\chi^\ell)-\ell(1-\chi^L)}{(1-\chi)(1-\chi^L)}\right)}.
\end{align}

We now distinguish between two different cases for $h_0$: either we
are free to choose the optimal $h_0\in\rset_+$, or we have an upper
bound on the coarsest mesh $h_0$. The first, idealized, situation will
allow us to obtain explicit expressions for the optimal splitting
parameter $\theta$ and the asymptotic work, and we start by
considering this case. We return to the other case at the end of this section.
\paragraph{Unconstrained optimization of $h_0$}
We take $h_1,\ldots,h_L$ given by~\eqref{eq:h_inner}
and~\eqref{eq:h_L} and set the derivative of~\eqref{eq:tmp_target}
with respect to $h_0$ equal to zero. This optimality condition becomes (after some straightforward simplifications)
\begin{align*}
  -\frac{d\gamma}{2}\frac{\sqrt{\vg}}{h_0^{1+d\gamma/2}}
  +\frac{q_2}{2}\sqrt{Q_S}\frac{h_0^{q_2/2-1}}{h_1^{d\gamma/2}}
  & = 0,
\end{align*}
which, since all parameters are positive, is equivalent to
\begin{align*}
  h_1 & =
  \left(
    \frac{\chi^2 Q_S}{\vg}
  \right)^{\frac{1}{d\gamma}}
  h_0^{1+\chi}.
\end{align*}
Combining this expression for $h_1$ with the one in~\eqref{eq:h_inner}
and solving for $h_0$ gives
\begin{align*}
  h_0 & =
  h_L^{\frac{1-\chi}{1-\chi^{L+1}}}
  \left(\frac{\vg}{Q_S}\right)^{\frac{1}{d\gamma}\frac{1-\chi^L}{1-\chi^{L+1}}}
  \chi^{-\frac{2}{d\gamma}\frac{1}{1-\chi}\left(L\frac{1-\chi}{1-\chi^{L+1}}-\chi\frac{1-\chi^L}{1-\chi^{L+1}}\right)},
\end{align*}
which after substituting back into~\eqref{eq:h_inner} and using \eqref{eq:h_L} yields~\eqref{eq:hl_optimal}.
Finally substituting these optimal mesh sizes into~\eqref{eq:opt_M_ell} yields~\eqref{eq:Ml_optimal}.

\paragraph{Optimal splitting parameter $\theta$}
Now the sequences $\{h_\ell\}_{\ell=0}^L$ and $\{M_\ell\}_{\ell=0}^L$
are determined in terms of the still not optimized $L$ and $\theta$ as
well as measurable model parameters.
The work per level in~\eqref{eq:model_work} becomes
\begin{align*}
  \frac{M_\ell}{h_\ell^{d\gamma}} & =
  \left(\frac{C_\alpha}{\theta\tol}\right)^2
  \left(
    \frac{Q_W}{(1-\theta)\,\tol}
  \right)^{\frac{1}{\eta}\frac{1-\chi}{1-\chi^{L+1}}}
   V_0
  \left(\frac{Q_S}{V_0}\right)^{\left\{\frac{1-\chi^L}{1-\chi^{L+1}}\right\}}
  \\
  &\quad\cdot
  \frac{1-\chi^{L+1}}{1-\chi}
      \chi^{\left\{-\frac{2\chi}{1-\chi}\frac{1-\chi^L}{1-\chi^{L+1}}
      +L\frac{1+\chi^{L+1}}{1-\chi^{L+1}}\right\}}
  \chi^{-\ell}.
\end{align*}
Since the only $\ell$-dependent factor in the right hand side is the
last one, $\chi^{-\ell}$, and using
$\sum_{\ell=0}^L\chi^{-\ell}=\chi^{-L}(1-\chi^{L+1})/(1-\chi)$,
the total work in~\eqref{eq:model_work} becomes
\begin{align}
  \label{eq:work_tmp_A}
  \work{(L,\theta,\tol)} & =
  w_1\left(L,\tol\right) \,
  w_2(L) \,
  f\left(L,\theta\right)
  \left(\frac{1-\chi^{L+1}}{1-\chi}\right)^2,
\end{align}
with
\begin{subequations}
  \label{eq:work_factors}
  \begin{align}
    \label{eq:w_1}
    w_1\left(L,\tol\right) & =
    \tol^{-\left(2+\frac{1}{\eta}\frac{1-\chi}{1-\chi^{L+1}}\right)}, \\
    \label{eq:w_2}
    w_2(L) & =
    C_\alpha^2 \,
    V_0
    \left(\frac{Q_S}{V_0}\right)^{\left\{\frac{1-\chi^L}{1-\chi^{L+1}}\right\}}
    Q_W^{\left\{\frac{1}{\eta}\frac{1-\chi}{1-\chi^{L+1}}\right\}}
    \\
    \nonumber
        &\quad\cdot
    \chi^{\left\{-\frac{2\chi}{1-\chi}\frac{1-\chi^L}{1-\chi^{L+1}}+2L\frac{\chi^{L+1}}{1-\chi^{L+1}}\right\}},
    \\
    \label{eq:def_f}
    f\left(L,\theta\right) & =
    \frac{1}{
      \theta^2
      \left(1-\theta\right)^{\frac{1}{\eta}\frac{1-\chi}{1-\chi^{L+1}}}
    }.
  \end{align}
\end{subequations}
Thus given the value of $L$ the dependence on the splitting parameter $\theta$ is
straightforward, and the minimal work for a given $L$ is obtained
with the minimizer of~\eqref{eq:def_f}, namely \eqref{eq:theta_opt}.
With this optimal splitting parameter $\theta$
in~\eqref{eq:work_tmp_A} the total work as a function of the yet
to be determined parameter $L$ and the tolerance is
\begin{align}
  \label{eq:work_theta_opt}
  \work{\left(L,\tol\right)} & =
  w_1(L,\tol) \, w_2(L) \, w_3 (L),
\end{align}
with
\begin{align}
  \label{eq:def_w_3}
  w_3(L) & =
  \left(\frac{1}{2\eta}\right)^2
  \left(
    1+2\eta\frac{1-\chi^{L+1}}{1-\chi}
  \right)^{2\left(1+\frac{1}{2\eta}\frac{1-\chi}{1-\chi^{L+1}}\right)}.
\end{align}

\paragraph{Optimal number of levels}
\label{s:app_optimal_L}

The optimal integer $L$ seems impossible to find
analytically.
In practical computations we instead perform an
extensive search over a small range of integer values.
In the analysis below we treat $L$ as a real parameter to obtain the
bounds~\eqref{eq:L_bound_TOL} that delimit the range of integer
values that must be tested, and allow a complexity analysis as
$\tol\to 0$ without an exactly determined $L$.

Treating $L$ as a real parameter, we differentiate the work~\eqref{eq:work_theta_opt} with respect to $L$ to obtain
\begin{align*}
    \frac{\partial \work}{\partial L} & =
  \frac{\partial w_1}{\partial L} \, w_2 \, w_3 +
  w_1 \, \frac{\partial w_2}{\partial L} \, w_3 +
  w_1 \, w_2 \, \frac{\partial w_3}{\partial L},
\end{align*}
where, introducing the shorthand
\begin{align}
  \label{eq:def_xi}
  \xi(L) & = 2\eta\frac{1-\chi^{L+1}}{1-\chi}
    && \text{for } L\in[0,\infty),
\end{align}
and using the constants $c_1$ and $c_2$ in \eqref{eq:L_bounds_consts}
we write
\begin{subequations}
  \label{eq:w_derivs}
  \begin{align}
  \frac{\partial w_1}{\partial L} & =
  w_1(L,\tol)\,\frac{\log{(\chi)}\,\chi^{L+1}}{1-\chi^{L+1}}
  \frac{2}{\xi(L)} \log{\left(\tol^{-1}\right)}, \\
  \frac{\partial w_2}{\partial L} & =
  w_2(L)\,
  \frac{\log{(\chi)}\,\chi^{L+1}}{1-\chi^{L+1}} \frac{2}{\xi(L)}
  \left(
    c_1
        + -c_2(L+1)
    + \xi(L)
  \right),\\
  \frac{\partial w_3}{\partial L} & =
  w_3(L)\,\frac{\log{(\chi)}\,\chi^{L+1}}{1-\chi^{L+1}} \frac{2}{\xi(L)}
  \big(
    \log{\left(1+\xi(L)\right)}-\xi(L)
  \big),
  \end{align}
\end{subequations}
so that
\begin{align}
  \label{eq:dWdL}
  \frac{\partial \work}{\partial L}(L,\tol) & =
  u(L,\tol)v(L,\tol),
\end{align}
with
\begin{align*}
    u(L,\tol) & =
  \work{(L,\tol)}\,
  \frac{\log{(\chi)}\,\chi^{L+1}}{1-\chi^{L+1}}
  \frac{2}{\xi(L)}, \\
    v(L,\tol) & =
    \log{\left(\tol^{-1}\right)} +
    c_1
    + -c_2 (L+1) +
    \log{\left(1+\xi(L)\right)}.
\end{align*}
Clearly $u(L,\tol)<0$ for all $\chi\in\rset_+\setminus\{1\}$ so the
sign of $\partial \work/\partial L$ is the opposite of the sign of
$v(L,\tol)$.
For a fixed $\chi\in\rset_+\setminus\{1\}$ we have
\begin{align*}
  v(L,\tol) > 0 & \Leftrightarrow
  L+1 < \frac{1}{c_2}
  \left(\log{\left(\tol^{-1}\right)}+c_1+\log{\left(1+\xi(L)\right)}\right),
\end{align*}
and, since $\xi(L)\geq\xi(0)=2\eta$,
\begin{align}
  \label{eq:low_bound_L_prel}
  L+1 < \frac{1}{c_2}
  \left(\log{\left(\tol^{-1}\right)}+c_1+\log{\left(1+2\eta\right)}\right)
  & \Rightarrow
  v(L,\tol) > 0 \Leftrightarrow \frac{\partial \work}{\partial L} < 0.
\end{align}
For the opposite inequality,
\begin{align*}
  v(L,\tol) < 0 & \Leftrightarrow
  L+1 > \frac{1}{c_2}
  \left(\log{\left(\tol^{-1}\right)}+c_1+\log{\left(1+\xi(L)\right)}\right),
\end{align*}
we distinguish between the cases $0<\chi<1$ and $1<\chi$. When
$0<\chi<1$ we have the upper bound $\xi(L)<\frac{2\eta}{1-\chi}$ and
consequently
\begin{align}
  \label{eq:upp_bound_L_prel1}
  L+1 > \frac{1}{c_2}
  \left(\log{\left(\tol^{-1}\right)}+c_1+\log{\left(1+\frac{2\eta}{1-\chi}\right)}\right)
  & \Rightarrow \frac{\partial \work}{\partial L} > 0, && \chi\in(0,1).
\end{align}
In contrast $\xi(L)$ is unbounded when $1<\chi$ but, since the
definitions of $\chi$ and $\eta$ and the relation between strong and
weak convergence orders implies that $2\eta\geq\chi$, we have
\begin{align*}
  \log{\left(1+\xi(L)\right)} & <
  \log{\left(\frac{2\eta}{\chi-1}\right)}
  + (L+1)\log{\chi},
\end{align*}
and
\begin{align*}
  c_2 & \geq \frac{\chi}{\chi-1}\log{\chi},
\end{align*}
which gives the bound
\begin{align*}
  \frac{1}{c_2} \log{\left(1+\xi(L)\right)} & <
  \frac{\chi-1}{\chi}(L+1) +
  \frac{1}{c_2} \log{\left(\frac{2\eta}{\chi-1}\right)}.
\end{align*}
Hence
\begin{align*}
  L+1 - \frac{1}{c_2} \log{\left(1+\xi(L)\right)} & >
  \frac{L+1}{\chi}
  - \frac{1}{c_2} \log{\left(\frac{2\eta}{\chi-1}\right)},
\end{align*}
and it follows that
\begin{align}
  \label{eq:upp_bound_L_prel2}
  L+1 > \frac{\chi}{c_2}
  \left(\log{\left(\tol^{-1}\right)}+c_1+\log{\left(\frac{2\eta}{\chi-1}\right)}\right)
  & \Rightarrow \frac{\partial \work}{\partial L} > 0, && \chi\in(1,\infty).
\end{align}

Combining~\eqref{eq:low_bound_L_prel} with \eqref{eq:upp_bound_L_prel1}
and \eqref{eq:upp_bound_L_prel2}, we obtain the bounds \eqref{eq:L_bound_TOL}.

\paragraph{Optimal hierarchies with an upper bound on $h_0$}
Practical computations will impose an upper limit on the mesh sizes,
$h_0\leq h_{\mathrm{max}}$. If the mesh sizes~\eqref{eq:hl_optimal}
violate such a bound, we must modify our analysis slightly. We now
consider $h_0$ given as one of the coarsest mesh sizes that can be
realized in the given discretization, and analyze the case $L\geq 1$.
Using the optimal mesh sizes~\eqref{eq:h_inner} yields
\begin{align*}
  \sqrt{\frac{h_{\ell-1}^{q_2}}{h_\ell^{d\gamma}}} & =
  h_0^{\frac{d\gamma}{2}\chi^L\frac{1-\chi}{1-\chi^L}}
  h_L^{-\frac{d\gamma}{2}\frac{1-\chi}{1-\chi^L}}
  \chi^{\left(\frac{L}{1-\chi^L}-\frac{\chi}{1-\chi}
      -\ell\right)},
\end{align*}

where the only $\ell$-dependent factor in the right hand side is the
last one, $\chi^{-\ell}$, so that the sum
in~\eqref{eq:work_factor_inner} is
\begin{align*}
  \sum_{\ell=1}^L \sqrt{\frac{h_{\ell-1}^{q_2}}{h_\ell^{d\gamma}}}
  & =
  \left(\frac{h_0^{\left(\chi^L\right)}}{h_L}\right)^{\frac{d\gamma}{2}\frac{1-\chi}{1-\chi^L}}
  \chi^{\left(\frac{L\chi^L}{1-\chi^L}-\frac{\chi}{1-\chi}\right)}
  \frac{1-\chi^L}{1-\chi}.
\end{align*}
In this sum only $h_L$ depends on $\theta$ through~\eqref{eq:h_L}.
Keeping $L$ fixed we wish to minimize the total
work, which by~\eqref{eq:total_work_app}--\eqref{eq:tmp_target} is
\begin{align*}
  \work(\hier) & =
  \left( \frac{C_\alpha}{ \theta\tol} \right)^2
  \left( \sqrt{\frac{\vg}{h_0^{d\gamma}}} +
    \sqrt{Q_S}
      \left(\frac{h_0^{\left(\chi^L\right)}}{h_L}\right)^{\frac{d\gamma}{2}\frac{1-\chi}{1-\chi^L}}
  \chi^{\left(\frac{L\chi^L}{1-\chi^L}-\frac{\chi}{1-\chi}\right)}
  \frac{1-\chi^L}{1-\chi}
  \right)^2,
\end{align*}
with respect to $\theta$. Letting
\begin{align*}
  \Delta & =
  \frac{1}{2\eta}\frac{1-\chi}{1-\chi^L},
\end{align*}
and
\begin{align*}
  C & =
  \sqrt{\frac{Q_S}{\vg}}
  h_0^{\frac{d\gamma}{2}\frac{1-\chi^{L+1}}{1-\chi^L}}
  \chi^{\left(\frac{L\chi^L}{1-\chi^L}-\frac{\chi}{1-\chi}\right)}
  \frac{1-\chi^L}{1-\chi}
  \left(\frac{Q_W}{\tol}\right)^\Delta,
\end{align*}
we obtain
\begin{align*}
  \work(\hier)
  & \propto
  \tilde{f}\left(\theta,L,h_0\right) =
  \frac{1}{\theta^2}\left(1+\frac{C}{\left(1-\theta\right)^\Delta}\right)^2,
\end{align*}
with the optimality condition
\begin{align*}
  \frac{\partial \tilde{f}}{\partial \theta}
  & =
  \frac{2}{\theta^2}\left(1+\frac{C}{\left(1-\theta\right)^\Delta}\right)
  \left(
    \frac{C\Delta}{\left(1-\theta\right)^{\Delta+1}}
      -\frac{1}{\theta}\left(1+\frac{C}{\left(1-\theta\right)^\Delta}\right)
  \right)
  = 0,
\end{align*}
where
\begin{align*}
  \frac{2}{\theta^2}\left(1+\frac{C}{\left(1-\theta\right)^\Delta}\right)
  & > 0.
\end{align*}

In this case when $h_0$ is constrained we no longer have an explicit
expression for the optimal $\theta$. However, using
\begin{align*}
  \frac{C\Delta}{\left(1-\theta\right)^{\Delta+1}}
  -\frac{1}{\theta}\left(1+\frac{C}{\left(1-\theta\right)^\Delta}\right)
  & <
  \frac{C}{\left(1-\theta\right)^\Delta}
  \left(\frac{\Delta}{1-\theta}-\frac{1}{\theta}\right),
\end{align*}
and that
\begin{align*}
  \frac{\Delta}{1-\theta}-\frac{1}{\theta} = 0
  & \Leftrightarrow
  \theta = \frac{1}{1+\Delta},
\end{align*}
we conclude that the optimal $\theta$ satisfies
\begin{align}
  \label{eq:low_bound_theta}
  \frac{1}{1+\Delta} & \leq \theta.
\end{align}
Similarly, from the inequality
\begin{align*}
  \frac{C\Delta}{\left(1-\theta\right)^{\Delta+1}}
  -\frac{1}{\theta}\left(1+\frac{C}{\left(1-\theta\right)^\Delta}\right)
  & >
  \frac{1}{\left(1-\theta\right)^\Delta}
  \left(\frac{C\Delta}{1-\theta}-\frac{1+C}{\theta}\right),
\end{align*}
and the relation
\begin{align*}
  \frac{C\Delta}{1-\theta}-\frac{1+C}{\theta} = 0
  & \Leftrightarrow
  \theta = \frac{1+C}{1+C+\Delta},
\end{align*}
we obtain an upper bound for $\theta$, namely
\begin{align}
  \label{eq:upp_bound_theta}
  \theta & \leq \frac{1+C}{1+C+C\Delta}.
\end{align}
Finally, combining~\eqref{eq:low_bound_theta} and~\eqref{eq:upp_bound_theta} we have the following bounds for the optimal $\theta$:
\begin{align}
  \label{eq:double_bound_theta}
  \left(1+\frac{1}{2\eta}\frac{1-\chi}{1-\chi^L}\right)^{-1}
  \leq \theta \leq
  \left(1+\frac{1}{2\eta}\frac{1-\chi}{1-\chi^L}\frac{C}{1+C}\right)^{-1},
\end{align}
where the upper bound has a non-trivial dependence on $\tol$ and $L$
through $C$.

\subsection{Heuristic optimization of geometric hierarchies}
\label{ss:app_geo}
This section motivates the results in Section~\ref{s:geo} and
Corollary~\ref{thm:opt_geo_work} where we optimized geometric
hierarchies defined by $h_\ell = h_0 \levsep^{-\ell}$ for given $h_0$
and $\levsep > 1$.  In this case, the work and variance models are in
\eqref{eq:work_var_geo_models} and $L$ is must satisfy the bias
constraint
\begin{align}
  L &\geq
  \frac{\log{(h_0)} - \frac{1}{q_1}\log{\left(\frac{(1-\theta)\tol}{Q_W}\right)}
     }{\log(\levsep)}.
    \label{eq:L_geo}
\end{align}
We distinguish between two cases:
\par{\textbullet\ $\chi=1$:}
Or equivalently $q_2 = d\gamma$. In this case, the total work defined
in \eqref{eq:total_work_app} simplifies to
\begin{equation}
\label{eq:work_geo_x1}
\work = \left( \frac{C_\alpha}{\theta \tol}\right)^2  \left(
\sqrt{V_0} h_0^\frac{
       -q_2}{2} + L \sqrt{Q_S} \beta^{\frac{q_2}{2}}  \right)^2,
\end{equation}
We make the simplification of treating $L$ as a real parameter and
substitute the lower bound of \eqref{eq:L_geo} in \eqref{eq:work_geo_x1} and
optimize with respect to $\levsep$ to get
$\levsep=\exp\left(\frac{2}{q_2}\right)$. Substituting this choice and \eqref{eq:L_geo_ceil}, the total
work satisfies
\[ \frac{\work}{\tol^{-2} \left( \log{\tol} \right)^2} \to \theta^{-2} C_\alpha^2 Q_S \exp(2) \left(\frac{1}{2\eta}\right)^2, \qquad \text{as } \tol \to 0.\]
Optimizing for $\theta$ suggests that $\theta \to 1$ as $\tol \to 0$ and \eqref{eq:opt_asymp_work_chi1} follows.

\par{\textbullet \ $\chi \neq 1$:} In this case, the total work defined in \eqref{eq:total_work_app} simplifies to
\begin{equation}
\label{eq:work_geo}
\work = \left( \frac{C_\alpha}{\theta \tol} \right)^2 h_0^{d\gamma(\chi-1)} \left(\sqrt{V_0} h_0^\frac{-q_2}{2} + \sqrt{Q_S} \frac{
        \left( 1-\levsep^\frac{L(d\gamma-q_2)}{2} \right)} {\levsep^{-\frac{d\gamma}{2}} - \levsep^{-\frac{q_2}{2}}}
        \right)^2,
    \end{equation}
for a given $L, h_0$ and $\theta$. Again, we make the simplification of treating $L$ as a real parameter and substitute the lower bound \eqref{eq:L_geo} to obtain
\[ \levsep^\frac{L(d\gamma-q_2)}{2} = \left( \frac{(1-\theta)\tol}{Q_W} \right)^\frac{\chi-1}{2\eta} h_0^\frac{d\gamma(1-\chi)}{2}, \]
for any $\levsep$. Substituting back in \eqref{eq:work_geo} and optimizing with respect
to $\levsep$ to minimize the work gives \eqref{eq:opt_beta}. Substituting this optimal $\levsep$ in \eqref{eq:work_geo} yields
\begin{equation}
\work = \left( \frac{C_\alpha}{\theta \tol} \right)^2 h_0^{d\gamma(\chi-1)} \left(\sqrt{V_0} h_0^\frac{-q_2}{2} + \sqrt{Q_S}
        \frac{\chi^\frac{\chi}{\chi-1}}{\chi-1} \left(1 - \chi^{-L}\right)
        \right)^2.
    \end{equation}
Asymptotically, using \eqref{eq:L_geo_ceil} as $\tol \to 0$ yields \eqref{eq:optimal_work} with the following constants
\begin{subequations}
\begin{align}
    C_1 &= (1-\theta)^\frac{\chi-1}{\eta}\theta^{-2} C_\alpha^2 Q_W^\frac{1-\chi}{\eta}  Q_S \left(
        \frac{\chi^\frac{\chi}{\chi-1}}{\chi-1} \right)^2, \\
    C_2 &=  \theta^{-2} C_\alpha^2 h_0^{d\gamma(\chi-1)} \left(\sqrt{V_0} h_0^\frac{-q_2}{2} + \sqrt{Q_S}
        \frac{\chi^\frac{\chi}{\chi-1}}{\chi-1} \right)^2.
\end{align}
\end{subequations}
Optimizing these constants with respect to $\theta$ yields
\eqref{eq:theta_opt_A} and substituting this and \eqref{eq:geo_opt_h0}
back yields \eqref{eq:optimal_work_C1} and \eqref{eq:C_2_geo} for
$C_1$ and $C_2$, respectively.  This, as Remark~\ref{rem:asymb_work}
mentions, shows that the asymptotic computational complexities of
optimal non-geometric and geometric hierarchies are the same.

\bibliographystyle{spmpsci}

\end{document}